\documentclass[12pt]{amsart}
\usepackage{graphicx}
\numberwithin{equation}{section}

\newcommand{\lab}{\label}

\newcommand{\ben}{\begin{enumerate}}
\newcommand{\een}{\end{enumerate}}

\newcommand{\bea}{\begin{eqnarray}}
\newcommand{\ba}{\begin{array}}
\newcommand{\bean}{\begin{eqnarray*}}
\newcommand{\ea}{\end{array}}
\newcommand{\eea}{\end{eqnarray}}
\newcommand{\eean}{\end{eqnarray*}}
\newcommand{\beq}{\begin{equation}}
\newcommand{\eeq}{\end{equation}}
\newcommand{\bthm}{\begin{thm}}
\newcommand{\ethm}{\end{thm}}
\newcommand{\blem}{\begin{lem}}
\newcommand{\elem}{\end{lem}}
\newcommand{\bprop}{\begin{prop}}
\newcommand{\eprop}{\end{prop}}
\newcommand{\bcor}{\begin{cor}}
\newcommand{\ecor}{\end{cor}}
\newcommand{\bdfn}{\begin{dfn}}
\newcommand{\edfn}{\end{dfn}}
\newcommand{\brem}{\begin{rem}}
\newcommand{\erem}{\end{rem}}
\newcommand{\bpf}{\begin{proof}}
\newcommand{\epf}{\end{proof}}
\newcommand{\bfact}{\begin{fact}}
\newcommand{\efact}{\end{fact}}

\newcommand{\nl}{\newline}

\alph{enumii} \roman{enumiii}

\newtheorem{thm}{Theorem}[section]
\newtheorem{prop}[thm]{Proposition}
\newtheorem{lem}[thm]{Lemma}

\newtheorem{cor}[thm]{Corollary}
\newtheorem{dfn}[thm]{Definition}
\newtheorem{rem}[thm]{Remark}
\newtheorem{fact}[thm]{Fact}
\newtheorem{ex}[thm]{Example}

\def\N{{\mathbb N}}                      \def\R{{\mathbb R}}
\def\C{{\mathbb C}}                  \def\oc{\hat \C}

\def\1{1\!\!1}
\def\and{\text{ and }}

\def\Comp{\text{Comp}}        \def\diam{\text{\rm {diam}}}

      \def\Exp{\text{{\rm Exp}}}

\def\H{\text{{\rm H}}}     \def\HD{\text{{\rm HD}}}   
         
\def\re{\text{{\rm Re}}}    \def\im{\text{{\rm Im}}}

         \def\P{\text{{\rm P}}}     \def\Id{\text{{\rm Id}}}

\def\L{{\mathcal L}}

\def\a{\alpha}                \def\b{\beta}             \def\d{\delta}
               \def\e{\varepsilon}          
\def\g{\gamma}                           \def\l{\lambda}
\def\La{\Lambda}              \def\om{\omega}           
\def\Sg{\Sigma}               \def\sg{\sigma}
               \def\th{\theta}           
\def\ka{\kappa}

\def\bi{\bigcap}              
\def\({\bigl(}                \def\){\bigr)}
\def\lt{\left}                \def\rt{\right}

\def\ld{\ldots}               \def\bd{\partial}         \def\^{\tilde}

\def\es{\emptyset}            \def\sms{\setminus}
\def\sbt{\subset}             \def\spt{\supset}

\def\sp{\medskip}                     \def\nl{\newline}
\def\ov{\overline}            

\def\ni{\noindent}

\def\om{\omega}

\def\re{\text{{\rm Re}}}

\begin{document}

\title[]
{\bf\large {\Large R}eal Analyticity of Hausdorff Dimension for Expanding 
Rational Semigroups}
\date{April 20, 2009. Published in Ergodic Theory  Dynam. Systems (2010), Vol. 30, No. 2, 601-633.}
\author[\sc Hiroki SUMI]{\sc Hiroki SUMI}
%
\author[\sc Mariusz URBA\'NSKI]{\sc Mariusz URBA\'NSKI}
%
%
\thanks{Research of the second author supported in part by the
NSF Grant DMS 0400481.}
\thanks{\ \newline 
\noindent Hiroki Sumi\newline 
Department of Mathematics,
Graduate School of Science,
Osaka University, 
1-1 Machikaneyama,
Toyonaka,
Osaka, 560-0043, 
Japan\newline 
E-mail: sumi@math.sci.osaka-u.ac.jp\newline 
Web: http://www.math.sci.osaka-u.ac.jp/$\sim$sumi/\newline
\ \newline 
Mariusz Urba\'nski\newline Department of Mathematics,
 University of North Texas, Denton, TX 76203-1430, USA\newline  
E-mail: urbanski@unt.edu\newline
Web: http://www.math.unt.edu/$\sim$urbanski/}

\keywords{Complex dynamical systems, rational semigroups, expanding semigroups,
Julia set, Hausdorff dimension, 
Bowen parameter , topological pressure, structural stability} 

\begin{abstract}
We consider the dynamics of expanding semigroups generated by 
finitely many rational maps on the Riemann sphere. 
We show that for an analytic family of such semigroups, 
the Bowen parameter function is real-analytic and plurisubharmonic. 
Combining this with a result obtained by the first author, 
we show that if for each semigroup of such an analytic family 
of expanding semigroups satisfies the open set condition, 
then the Hausdorff dimension of the Julia set 
is a real-analytic and plurisubharmonic function of the parameter. Moreover, we
provide 
an extensive collection of examples of analytic families 
of semigroups satisfying all the above conditions and we analyze 
in detail the corresponding Bowen's parameters and Hausdorff dimension 
function.    
\end{abstract}
\maketitle
 Mathematics Subject Classification (2001). Primary 37F35; 
Secondary 37F15.

\section{Introduction}
\label{Introduction}
A {\bf rational semigroup} 
is a semigroup generated by a family of 
non-constant rational maps $g:\oc \rightarrow \oc $, 
where $\oc $ denotes the Riemann sphere, 
with the semigroup operation being functional composition. 
A polynomial semigroup is a semigroup generated by a 
family of non-constant polynomial maps on $\oc .$ 
Research on the dynamics of rational semigroups was initiated 
by A. Hinkkanen and G. J. Martin (\cite{HM, HM2}), 
who were interested in the role of the dynamics of polynomial semigroups 
while studying various one-complex-dimensional moduli spaces for discrete 
groups of M\"{o}bius transformations, and by F. Ren's group 
(\cite{ZR,GR}), who studied such semigroups from the perspective 
of random dynamical systems. 
  
 The theory of the dynamics of rational semigroups on $\oc $ 
has developed in many directions since the 1990s (\cite{HM, ZR, HM2,St1, St2, St3, SSS, sumihyp1, sumihyp2, 
hiroki1, sumi1, sumicorrection, sumi2, sumi06, sumirandom, sumikokyuroku, su0, su1, sumid1, sumid3, sumi07, SS, sumiintcoh, 
sumiprepare}). 
Since the Julia set $J(G)$ of a rational semigroup 
generated by finitely many elements $f_{1},\ldots ,f_{s}$ 
has {\bf backward self-similarity} i.e. 
\begin{equation}
\label{bsseq}
 J(G)=f_{1}^{-1}(J(G))\cup \cdots \cup f_{s}^{-1}(J(G)),
\end{equation}  
(See \cite{sumihyp1, hiroki1}), 
it can be viewed significant generalization and extension of 
both the theory of iteration of rational maps (see \cite{M}) 
and conformal 
iterated function systems (see \cite{MU}). The theory of the dynamics of 
rational semigroups borrows and develops tools 
from both of these theories. It has also developed its own 
unique methods, notably the skew product approach 
(see \cite{hiroki1, sumi1, sumi2, sumi06, sumid1, sumid3, sumi07, sumiprepare}). 
We remark that by (\ref{bsseq}), 
the analysis of the Julia sets of rational semigroups somewhat
resembles 
 ``backward iterated functions systems'', however since each map 
$f_{j}$ is not in general injective (critical points), some 
qualitatively different extra effort in the cases of semigroups is needed.

The theory of the dynamics of rational semigroups is intimately 
related to that of the random dynamics of rational maps. 
For the study of random complex dynamics, the reader may 
consult \cite{FS,Bu1,Bu2,BBR,Br,GQL}. The deep relation between these fields 
(rational semigroups, random complex dynamics, and (backward) IFS) 
is explained in detail in the subsequent papers
 (\cite{sumirandom, sumikokyuroku, sumid1, sumid3, sumi07, sumiintcoh, sumiprepare}) of the first author.

 In this paper, we analyze in detail the Hausdorff dimension 
 of Julia sets of expanding rational semigroups. 
 Our approach utilize the powerful tool of thermodynamic formalism, 
 developed in \cite{sumi2} and a version of Bowen's formula 
 for the Hausdorff dimension of Julia sets, also proved 
 in \cite{sumi2}. We introduce Bowen's parameter as the unique 
 zero of the pressure function. This is an invariant 
 of the generator systems of the semigroup. We then develop 
 a finer analysis of holomorphic families of Perron-Frobenius 
 type operators, and eventually apply Kato-Rellich perturbation 
 theory (\cite{Ka}) to get real-analyticity of the pressure function, 
 as depending on the complex parameter. Then the Implicit  
 Function Theorem completes the task. Bowen's formula, which is 
 mentioned above, identifies the Hausdorff dimension of the Julia set 
 with Bowen's parameter, whenever in addition the open set condition 
 is satisfied. We thus obtain that under these assumptions the 
 Hausdorff dimension function depends in a real-analytic manner 
 on the parameter. We also show that 
 Bowen's parameter function is real-analytic and plurisubharmonic, 
 even if we do not assume the open set condition. The real analyticity
 of Hausdorff dimension, or, in a sense, more accurately, of Bowen's
 parameter, goes back to
 Ruelle's paper \cite{ruelle}, where hyperbolic rational functions were
 considered. The reader may also consult \cite{zinsmeister} and
 \cite{zdunik}. Our approach stems from that in the \cite{zdunik}
 paper. We develop it and work out techniques to deal with a
 qualitatively different case of semigroups. 
  
  Our article ends with a collection of examples illustrating variety 
  of behavior of the Hausdorff dimension function, 
  Bowen's parameter function (ex. it can be strictly less than $2$ or 
  strictly bigger than $2$ on an open set of multi-maps), expandingness 
and the open set condition.     
 
 We remark that as illustrated in \cite{sumikokyuroku,sumiprepare}, 
 estimating the Hausdorff dimension of the Julia sets of 
 rational semigroups plays an important role when we 
 investigate random complex dynamics and its associated 
 Markov process on $\oc .$ For example, 
 when we consider the random dynamics of a compact 
 family $\Gamma $ of polynomials 
 of degree greater than or equal to two,  
 then the function $T_{\infty }:\oc \rightarrow 
 [0,1]$ of probability of tending to $\infty \in \oc $ 
 varies only inside the Julia set of rational semigroup  
 generated by $\Gamma $, and under some condition, 
 this $T_{\infty }:\oc \rightarrow [0,1]$ is continuous in $\oc .$
If the Hausdorff dimension of the Julia set is strictly less than two, 
then it means that $T_{\infty }:\oc \rightarrow [0,1]$ is a 
complex version of devil's staircase (Cantor function) 
(\cite{sumirandom,sumikokyuroku,sumiprepare}). 
\section{Preliminaries and the main results}
In this section we introduce the notation and basic definitions. 
We also formulate our main results. 
Their proofs will be concluded in Section~\ref{Real}. 

Throughout the paper, we frequently follow the notation 
from \cite{hiroki1} and \cite{sumi2}. 
\begin{dfn}[\cite{HM,ZR,GR,HM2}] 
A ``rational semigroup" $G$ is a semigroup generated by a family of 
non-constant 
rational maps $g:\oc \rightarrow \oc$,\ where $\oc $ denotes the 
Riemann sphere,\ with the semigroup operation being functional 
composition. 
A ``polynomial semigroup'' is a semigroup generated by a 
family of non-constant polynomial maps on $\oc .$ 
For a rational semigroup $G$,\ we set 
$F(G):=\{ z\in \oc \mid G \mbox{ is normal in a neighborhood of } z\} 
$ and $J(G):=\oc \setminus F(G).$ 
$F(G)$ is called the {\bf Fatou set} of $G$ and $J(G)$ is called the 
{\bf Julia set} of $G.$ 
If $G$ is generated by a family $\{ f_{i}\} _{i}$,\ 
then we write $G=\langle f_{1},f_{2},\cdots \rangle .$ 
\end{dfn} 
For the study of the dynamics of rational semigroups, 
see \cite{HM,ZR,GR,HM2,St1,St2,St3,SSS, sumihyp1,sumihyp2,
hiroki1,sumi1,sumicorrection,sumi2,sumi06, sumirandom, sumikokyuroku, su0, su1, 
sumid1, sumid3, sumi07,SS,sumiintcoh, sumiprepare}, etc. 

\begin{dfn}
For each $s\in \Bbb{N}$, 
let $\Sigma _{s}:=\{ 1,\ldots ,s\} ^{\Bbb{N}}$ be the 
space of one-sided sequences of $s$-symbols endowed with the 
product topology. This is a compact metric space. 
We denote by {\em Rat} the set of all non-constant 
rational maps on $\oc $ endowed with 
the topology induced by uniform convergence on $\oc .$ 
Note that {\em Rat} has countably many connected 
components. In addition, each connected component 
$U$ of {\em Rat} is an open subset of {\em Rat} and 
$U$ has a structure of a finite dimensional complex manifold.  
For each $f=(f_{1},\ldots ,f_{s})\in (\mbox{{\em Rat}})^{s}$, 
we define a map   
$$
\tilde{f}:\Sg_{s}\times \oc \rightarrow \Sg_{s}\times \oc 
$$ 
by the formula 
$$
\tilde{f}(\om,z)=(\sg (\om ),\ f_{\om_{1}}(z)),
$$
where $(\om,z)\in \Sg _{s}\times \oc,\ \om=(\om_{1},\om_{2},\ldots ),$ and 
$\sg :\Sigma _{s}\rightarrow \Sg _{s}$ denotes the shift map. 
This $\tilde{f} :\Sigma _{s}\times \oc \rightarrow 
\Sigma _{s}\times \oc $ is called the {\bf skew product map} associated 
with the multi-map $f=(f_{1},\ldots ,f_{s})\in (\mbox{{\em Rat}})^{s} .$  
We denote by $\pi _{1}:\Sigma _{s}\times \oc \rightarrow \Sigma _{s}$ 
the projection onto $\Sg_{s}$ and $\pi_{2}:\Sg _{s}\times\oc\rightarrow\oc$ 
the projection onto $\oc $. That is, $\pi _{1}(\om ,z)=\om $ and 
$\pi _{2}(\om ,z)=z.$  For each $n\in \Bbb{N} $ and $(\om ,z)\in 
\Sigma _{s}\times \oc $, we put 
$$
(\tilde{f}^{n})'(\om ,z):= (f_{\om_{n}}\circ \cdots \circ f_{\om _{1}})'(z).
$$ 
We put $J_{\om }(\tilde{f}):=\{ z\in \oc \mid 
 \{ f_{\om _{n}}\circ \cdots \circ f_{\om _{1}}\} _{n\in \Bbb{N}} \mbox{ is 
 not normal in each neighborhood of } z\} $ for each $\om \in 
 \Sigma _{s}$ and we set 
$$
J(\tilde{f}):= \overline{\cup _{w\in \Sigma _{s}}\{ \om \} 
\times J_{\om }(\tilde{f}) },
$$ 
where the closure is taken in the product space 
$\Sigma _{s}\times \oc .$ This is called the 
{\bf Julia set} of the skew product map $\tilde{f}.$  
Moreover, we set 
$F(\tilde{f}):=(\Sigma _{s}\times \oc )\setminus J(\tilde{f}).$ 
 Furthermore, we set  
$
\deg(\tilde{f}):=\sum _{j=1}^{s}\deg(f_{j}).
$
\end{dfn}

\begin{rem}
\label{rem1} 
By definition,\ 
 $J(\tilde{f})$ 
 is compact. Furthermore,\ if we set 
 $G=\langle f_{1},\ldots ,f_{s}\rangle $, then 
 by  \cite[Proposition 3.2]{hiroki1},\ we have all of the following:  
\begin{enumerate}
\item  
 $J(\tilde{f})$ is completely invariant under $\tilde{f}$;\ 
\item 
 $\tilde{f}$ is an open map on $J(\tilde{f})$;\ 
\item 
if 
 $\sharp J(G)\geq 3$ and $E(G):= \{ z\in \oc \mid 
 \sharp \cup _{g\in G}g^{-1}\{ z\}<\infty \} $ is included in 
 $F(G)$,  then  $(\tilde{f},J(\tilde{f}))$ is topologically exact; 
\item  
 $J(\tilde{f})$ is equal to the closure of 
 the set of repelling periodic points of 
 $\tilde{f}$ if 
 $\sharp J(G)\geq 3$,\ where we say that a periodic point 
 $(\om ,z)$ of $\tilde{f}$ with 
 period $n$ is repelling if $|(\tilde{f}^{n})'(\om ,z)|>1$; and       
\item $\pi _{2}(J(\tilde{f}))=J(G).$ 
\end{enumerate}
\end{rem} 
\begin{dfn}[\cite{sumi2}]
A finitely generated 
rational semigroup $G=\langle f_{1},\ldots ,f_{s}\rangle $ 
is said to be expanding provided that 
$J(G)\ne \es $ and the skew product map 
$\tilde{f}:\Sg_{s} \times \oc \rightarrow \Sg _{s}\times \oc $ 
associated with 
$f=(f_{1},\ldots ,f_{s}) $ is expanding along 
fibers of the Julia set $J(\tilde{f})$, 
meaning that there exists $\eta >1$ and 
$C\in(0,1]$ such that for all $n\ge 1$,
\begin{equation}
\label{1112505}
\inf \{ \| (\tilde{f}^{n})'(z)\|: z\in J(\tilde{f})\} \ge C\eta ^{n},  
\end{equation}  
where we mean in the above formula $\|\cdot \| $ to denote 
the absolute value of the spherical derivative. 
\end{dfn} 
\begin{dfn}
Let $G$ be a rational semigroup. 
We set 
$$P(G):=\overline{\cup _{g\in G}\{ \mbox{all critical values of }
g:\oc \rightarrow \oc \}} \ (\subset \oc )$$ 
and this is called the {\bf postcritical set} of $G$. 
A rational semigroup $G$ is said to be {\bf hyperbolic} if 
$P(G)\subset F(G).$ 
\end{dfn}

\begin{rem}
\label{exphyplem}
Let $G=\langle f_{1},\ldots ,f_{s}\rangle $ be a rational semigroup 
such that 
there exists an element $g\in G$ with $\deg (g)\geq 2$ and 
such that each M\"{o}bius transformation in $G$ is loxodromic. 
Then,  
it was proved in \cite{sumihyp2} that 
$G$ is expanding if and only if $G$ is hyperbolic. 
\end{rem}
\begin{dfn}
We define
$$\Exp (s):=\{ (f_{1},\ldots ,f_{s})\in (\mbox{{\em Rat}})^{s}
\mid \langle f_{1},\ldots ,f_{s}\rangle \mbox{ is expanding} \}.
$$  
We also set $\Sigma _{s}^{\ast }:= \cup _{j=1}^{\infty }
\{ 1,\ldots ,s\} ^{j}$ (disjoint union). 
For every 
$\om\in\Sg_s\cup\Sg_s^*$ let $|\om|$ be the length of $\om .$
For each $f=(f_{1},\ldots ,f_{s})\in (\mbox{{\em Rat}})^{s}$ 
and each $\om =(\om _{1},\ldots ,\om _{n})\in \Sigma _{s}^{\ast }$,  
we put $f_{\om }:= f_{\om _{n}}\circ \cdots \circ f_{\om _{1}}.$ 
\end{dfn}
Then we have the following.
\begin{lem}
\label{expopenlem}
$\Exp(s)$ is an open subset of {\em (Rat}$)^{s}.$ 
\end{lem}
{\sl Proof.}
Let $f=(f_{1},\ldots ,f_{s})\in \Exp (s).$ 
Then, by (\ref{1112505}) and the fact 
$\pi _{2}(J(\tilde{f}))=J(\langle f_{1},\ldots ,f_{s}\rangle )$ (Remark~\ref{rem1}), 
there exists an $n\in \N $ such that 
\begin{equation}
\label{expopenlemeq1}
\inf \{ \| (f_{\om })'(y)\| :\  
\om \in \Sigma _{s}^{\ast }, |\om |=n, 
y\in f_{\om }^{-1}(J(\langle f_{1},\ldots ,f_{s}\rangle ))\} \geq 3.
\end{equation} 
For each subset $A$ of $\oc $ and $r>0$, we denote by   
$B(A,r)$ the $r$-neighborhood of $A$ with respect to 
the spherical distance on $\oc .$ Let $\epsilon >0$ be any small number. 
Then, by (\ref{expopenlemeq1}), for each $\om \in \Sigma _{s}^{\ast }$ 
with $|\om |=n$, 
\begin{equation}
f_{\om }^{-1}(B(J(\langle f_{1},\ldots ,f_{s}\rangle ), \epsilon ))
\subset B(J(\langle f_{1},\ldots ,f_{s}\rangle ), \epsilon /2). 
\end{equation} 
Hence, there exists a neighborhood $U$ of $f$ in 
$(\mbox{Rat})^{s}$ such that 
for each $g=(g_{1},\ldots ,g_{s})\in U$ and each $\om $ with $|\om |=n$, 
\begin{equation}
\label{expopenlemeq1.5}
g_{\om }^{-1}(B(J(\langle f_{1},\ldots ,f_{s}\rangle ), \epsilon ))
\subset B(J(\langle f_{1},\ldots ,f_{s}\rangle ),\epsilon )
\end{equation} 
 and 
\begin{equation}
\label{expopenlemeq2}
\| (g_{\om })'(y)\| > 3/2 \mbox{ for each } 
y\in g_{\om }^{-1}(B(J(\langle f_{1},\ldots ,f_{s}\rangle ),\epsilon )).  
\end{equation}
By (\ref{expopenlemeq1.5}), 
for each $g\in U$, 
$\langle g_{1},\ldots ,g_{s}\rangle $ is normal in 
$\oc \setminus \overline{B}(J(\langle f_{1},\ldots ,f_{s}\rangle ),\epsilon ).$ 
Hence,  
it follows that 
for each $g\in U$, 
\begin{equation}
\label{expopenlemeq3}
J(\langle g_{1},\ldots ,g_{s}\rangle )\subset 
\overline{B}(J(\langle f_{1},\ldots ,f_{s}\rangle ), \epsilon ).
\end{equation} 
By (\ref{expopenlemeq2}), (\ref{expopenlemeq3}), and 
the fact $\pi _{2}(J(\tilde{g}))=J(\langle g_{1},\ldots ,g_{s}\rangle )$, 
we obtain that for each $g\in U$, 
$\langle g_{1},\ldots ,g_{s}\rangle $ is expanding. 
We are done. 
\qed 

\begin{dfn}
Let $f=(f_{1},\ldots ,f_{s})\in \Exp(s)$ and 
let $\tilde{f}:\Sigma _{s}\times \oc \rightarrow 
\Sigma _{s}\times \oc $ be the skew product map 
associated with $f=(f_{1},\ldots ,f_{s}) .$ 
For each $t\in \R $, 
let $P(t,f)$ be the topological pressure of 
the potential $\varphi _(z):= -t\log \| \tilde{f}'(z)\|$ 
with respect to the map 
$\tilde{f}:J(\tilde{f})\rightarrow J(\tilde{f}).$  
(For the definition of the topological pressure, 
see \cite{pubook}.) 
We denote by $\delta (f)$ the unique zero 
of $t\mapsto P(t,f).$ (Note that the existence and 
the uniqueness of the zero of $P(t,f)$ was shown in 
\cite{sumi2}.) This $\delta (f)$ is called the 
{\bf Bowen parameter} of $f=(f_{1},\ldots ,f_{s})\in 
\Exp(s).$  
\end{dfn}

In this paper, we consider the following situation: 
\begin{dfn}
Let $\Lambda $ be a finite dimensional complex manifold. 
For each $j=1,\ldots ,s$ and each $\lambda \in \Lambda $, 
suppose that there exists a rational map 
$f_{\lambda ,j}:\oc \rightarrow \oc .$  
For each $\lambda \in \Lambda $, we set 
$G_{\lambda }:= \langle f_{\lambda ,1},\ldots ,f_{\lambda ,s}\rangle .$ 
The collection $\{ G_{\lambda }\} _{\lambda \in \Lambda }$ is 
called an analytic family of rational semigroups 
if 
\begin{itemize}
\item[(a)] 
For every $1\le j\le s$ and every $z\in\C$,  
$(z,\l) \mapsto f_{\l ,j}(z)$ is a holomorphic map from 
$\oc \times \Lambda $ to $\oc .$ 
\end{itemize}
Furthermore, the collection  
$\{ G_{\lambda }\} _{\lambda \in \Lambda }$ 
is called an analytic family of expanding rational semigroups if 
$\{ G_{\lambda }\} _{\lambda \in \Lambda }$ is an analytic family 
of rational semigroups and   
for all $\l \in \Lambda $, 
$G_{\l }$ is expanding. 
\end{dfn}
\begin{ex}
Let $\Lambda $ be a connected component of $\Exp(s).$ For each 
$f=(f_{1},\ldots ,f_{s})\in \Lambda $, let 
$G_{f}:=\langle f_{1}, \ldots ,f_{s}\rangle .$ Then, 
by Lemma~\ref{expopenlem}, 
$\{G_{f}\} _{f\in \Lambda }$ is an analytic family of 
expanding rational semigroups.  
\end{ex}

We will give a large collection of examples of analytic families of expanding 
rational semigroups in Section~\ref{Examples}.
 
In order to state the main results, we need the following notation.
\begin{dfn}[\cite{K}]
Let $X$ be a finite dimensional complex manifold. 
An upper semicontinuous function $u:X\rightarrow \Bbb{R}\cup \{ -\infty \} $ 
is said to be plurisubharmonic if for each holomorphic 
map $\varphi :\Bbb{D}\rightarrow  X$, 
where $\Bbb{D}:=\{ z\in \Bbb{C}:|z|<1\} $, 
the function $u\circ \varphi :\Bbb{D}\rightarrow \Bbb{R}\cup \{ -\infty \} $ 
is subharmonic. 
A function $v:X\rightarrow \Bbb{R}\cup \{ +\infty \} $ is 
said to be plurisuperharmonic if 
$-v :X\rightarrow \Bbb{R}\cup \{ -\infty \} $ is 
plurisubharmonic.    
\end{dfn} 
\begin{dfn}
For any subset $A$ of $\oc $, 
we denote by $\HD(A)$ the Hausdorff dimension of $A$ with respect to 
the spherical distance on $\oc .$ 
\end{dfn}
\begin{dfn}
Let $f=(f_{1},\ldots ,f_{s})\in (\mbox{{\em Rat}})^{s}$ be an element 
and let $G=\langle f_{1},\ldots ,f_{s}\rangle $.  
Let also $U$ be a non-empty open set in $\oc .$ 
We say that $f$ (or $G$) satisfies the open set condition (with $U$) 
if $\cup _{j=1}^{s}f_{j}^{-1}(U)\subset U$ and 
$f_{i}^{-1}(U)\cap f_{j}^{-1}(U)=\emptyset $ for 
each $(i,j)$ with $i\neq j.$ 
There is a stronger condition. Namely, we say that 
$f$ (or $G$) satisfies the separating open set condition 
(with $U$) 
if $\cup _{j=1}^{s}f_{j}^{-1}(U)\subset U$ and 
$f_{i}^{-1}(\overline{U})\cap f_{j}^{-1}(\overline{U})=\emptyset $ for 
each $(i,j)$ with $i\neq j.$ 
Furthermore, we say that 
$f$ (or $G$) satisfies the 
strongly separating open set condition 
(with $U$) 
if $\cup _{j=1}^{s}f_{j}^{-1}(\overline{U})\subset U$ and 
$f_{i}^{-1}(\overline{U})\cap f_{j}^{-1}(\overline{U})=\emptyset $ for 
each $(i,j)$ with $i\neq j.$ 

\end{dfn}

The main purpose of this paper is to prove the following results:
\begin{thm}
\label{mainthmA}
({\bf Theorem A}) 
Let $\Lambda $ be a finite dimensional complex manifold. 
Let $\{ G_{\lambda }\} _{\lambda \in \Lambda }$ 
be an analytic family of expanding rational semigroups, 
where $G_{\l }=\langle f_{\l ,1},\ldots ,f_{\l ,s}\rangle .$   
For each $\l \in \Lambda $, 
we set $f_{\l }:= (f_{\l ,1},\ldots ,f_{\l ,s})\in (\mbox{{\em Rat}})^{s}.$ 
Then, the Bowen parameter function 
$\lambda \mapsto \delta (f_{\l })$ defined for all 
$\lambda \in \Lambda $, is real-analytic.   
Also, $(\lambda ,t)\mapsto P(t,f_{\lambda }), 
(\lambda ,t)\in \Lambda \times \Bbb{R},$ is real-analytic, 
$\lambda \mapsto 1/\delta (f_{\l }), \l \in \Lambda $, 
is plurisuperharmonic, $\l \mapsto \delta (f_{\l }), \l \in \Lambda $, 
is plurisubharmonic, and $\l \mapsto \log \delta (f_{\l }), 
\l \in \Lambda $, is plurisubharmonic. Furthermore, 
for a fixed $t\in \Bbb{R}$, the function 
$\lambda \mapsto P(t,f_{\lambda }),\lambda \in \Lambda ,$ is 
plurisubharmonic.   
\end{thm}

\begin{thm}
\label{mainthmB}
({\bf Theorem B})
Let $\Lambda $ be a finite dimensional complex manifold. 
Let $\{ G_{\lambda }\} _{\lambda \in \Lambda }$ 
be an analytic family of expanding rational semigroups, 
where $G_{\l }=\langle 
f_{\lambda ,1},\ldots ,f_{\lambda ,s}\rangle .$   
Suppose that for each $\lambda \in \Lambda $, 
$G_{\lambda }$ satisfies the open set condition i.e., 
for each $\lambda \in \Lambda $ there exists a non-empty 
open set $U_{\lambda }$ in $\oc $ such that 
$\cup _{j=1}^{s}f_{\lambda ,j}^{-1}(U_{\lambda })
\subset U_{\lambda }$ and 
$f_{\lambda ,i}^{-1}(U_{\lambda })\cap 
f_{\lambda ,j}^{-1}(U_{\lambda })=\emptyset $ for each 
$(i,j)$ with $i\neq j.$  
Then, 
the Hausdorff dimension function 
$\lambda \mapsto \HD(J(G_{\lambda })), \lambda \in \Lambda $, 
is real-analytic. Besides, 
$\l \mapsto 1/\HD(J(G_{\l })), \l \in \Lambda $, 
is plurisuperharmonic, 
$\l \mapsto \HD(J(G_{\l })), \l \in \Lambda $, 
is plurisubharmonic, 
and $\l \mapsto \log \HD(J(G_{\l })), \l \in \Lambda $, 
is plurisubharmonic. 
\end{thm}
\begin{rem}
There exist a number of elements $g=(g_{1},\ldots ,g_{s})\in 
\Exp(s)$ such that the Hausdorff dimension function 
$f=(f_{1},\ldots ,f_{s})\mapsto 
\HD(J(\langle f_{1},\ldots ,f_{s}\rangle ))$ 
is not continuous at $g.$ For example, 
let $g=(z^{2},z^{2})\in \Exp(2)$  
and for each $\l $ with $0<\l \leq 1$, let 
$f_{\l }:= (z^{2}, \l z^{2}) $ 
and $G_{\l }:=\langle z^{2}, \l z^{2}\rangle .$ 
Then, for each $0<\l <1$, 
$J(G_{\l })=\{ z\in \Bbb{C} : 1\leq |z|\leq 1/\l \} $ which implies 
$\HD(J(G_{\l }))=2.$ However, 
$J(G_{1})=J(\langle z^{2},z^{2}\rangle )=\{ z\in \Bbb{C} : |z|=1\} $ and 
$\HD(J(G_{1}))=1.$ 
\end{rem}
\begin{rem}
\label{oscrem} \ 
\begin{itemize}
\item  
There is a rich collection of finitely generated rational semigroups 
$G=\langle f_{1},\ldots ,f_{s}\rangle $ such that 
{\em int}$J(G)\neq \emptyset $ and 
$J(G)\neq \oc .$ For example, 
for each $\l \in \Bbb{C}$ with $0<|\l |<1$, 
$J(\langle z^{2}, \l z^{2}\rangle )=\{ z\in \Bbb{C}: 1\leq |z|\leq 1/|\l |\} .$ 
\item 
If a finitely generated 
rational semigroup $G=\langle f_{1},\ldots ,f_{s}\rangle $ 
satisfies the open set condition with 
an open set $U$, then by \cite[Corollary 3.2]{HM}, 
$J(G)\subset \overline{U}.$ 
If, in addition to the above, $J(G)\neq \overline{U}$, 
then {\em int}$J(G)=\emptyset $ 
(See \cite[Proposition 4.3]{sumi1}). 
\item Suppose a finitely generated 
rational semigroup 
$G=\langle f_{1},\ldots ,f_{s}\rangle $ satisfies 
the separating open set condition with the set $U.$   
Then, by \cite[Corollary 3.2]{HM},   
\cite[Theorem 2.3]{sumihyp2}, and 
\cite[Lemma 2.4]{hiroki1}, 
we have that {\em int}$J(G)=\emptyset $ and 
$J(G)$ is disconnected. In particular, 
$J(G)$ is a proper disconnected subset of $\overline{U}.$ 
\item 
If a finitely generated expanding rational semigroup 
$G=\langle f_{1},\ldots ,f_{s}\rangle $ 
satisfies the open set condition with 
the set $U$ and $J(G)\neq \overline{U}$, 
then by \cite[Theorem 1.25]{sumi06} and its proof, 
the Julia set $J(G)$ is porous and $\HD(J(G))<2.$ In particular, 
if an expanding rational semigroup 
$G=\langle f_{1},\ldots ,f_{s}\rangle $ satisfies the 
separating open set condition with the set $U$, then 
$J(G)$ is porous and $\HD(J(G))<2.$ 
\end{itemize} 
\end{rem}

The proofs of Theorem A and Theorem B are given 
in Section~\ref{Real}. They make use of the 
thermodynamic formalisms and 
the perturbation theory for bounded linear operators on 
Banach spaces. 

We give some additional remarks.
\begin{dfn}[\cite{sumi2}]
Let $G$ be a countable rational semigroup. 
For any $t\geq 0$ and $z\in \oc $, we 
set $S_{G}(z,t):=\sum _{g\in G}\sum _{g(y)=z}\| g'(y)\| ^{-t}$, 
counting multiplicities.  
We also set 
$S_{G}(z):= \inf \{ t\geq 0: S_{G}(z,t)<\infty \} $ 
(if no $t$ exists with $S_{G}(z,t) <\infty $, then we set 
$S_{G}(z):=\infty $). Furthermore, 
we set $s_{0}(G):= \inf \{ S_{G}(z): z\in \oc \} .$  
This $s_{0}(G)$ is called the {\bf critical exponent of the 
Poincar\'{e} series} of $G.$ 
\end{dfn}
\begin{dfn}[\cite{sumi2}]
Let $f=(f_{1},\ldots ,f_{s})\in (\mbox{{\em Rat}})^{s}$, $t\geq 0$, and $z\in \oc .$
 We put
$T_{f}(z,t):=\sum _{\om \in \Sigma _{s}^{\ast}}
\sum _{f_{\om }(y)=z}\| f_{\om }'(y)\| ^{-t}$, 
counting multiplicities.   
Moreover, we set 
$T_{f}(z):=\inf \{ t\geq 0:T_{f}(z,t)<\infty \} $ 
(if no $t$ exists with $T_{f}(z,t)<\infty $, then we set 
$T_{f}(z)=\infty $). 
Furthermore, we set 
$t_{0}(f):= \inf \{ T_{f}(z): z\in \oc \} .$  
This $t_{0}(f)$ is called the {\bf critical exponent of 
the Poincar\'{e} series} of $f=(f_{1},\ldots ,f_{s})\in 
(\mbox{{\em Rat}})^{s}.$ 
\end{dfn}
\begin{rem}
\label{strem}
Let $f=(f_{1},\ldots ,f_{s})\in (\mbox{{\em Rat}})^{s}$, $t\geq 0$ , 
$z\in \oc $ and let 
$G=\langle f_{1},\ldots ,f_{s}\rangle .$ 
Then, 
$S_{G}(t,z)\leq T_{f}(t,z), S_{G}(z)\leq T_{f}(z),$ and 
$s_{0}(G)\leq t_{0}(f).$ Note that 
 for almost every $f\in (\mbox{{\em Rat}})^{s}$ with 
 respect to the Lebesgue measure, 
 $G=\langle f_{1},\ldots ,f_{s}\rangle $ is a free 
 semigroup and so we have 
 $S_{G}(t,z)=T_{f}(t,z), S_{G}(z)=T_{f}(z), $ and 
$s_{0}(G)=t_{0}(f).$  
\end{rem}

\begin{dfn}
Let $G$ be a rational semigroup. 
Then, we set 
$$A(G):= \overline{\cup _{g\in G}g(\{ z\in \oc : \exists u\in G, 
u(z)=z, |u'(z)|<1\} )}.$$
\end{dfn}
\begin{lem}
\label{deltat0lem}
Let $f=(f_{1},\ldots ,f_{s})\in \Exp(s).$ 
Then $\delta (f)=t_{0}(f).$ 

{\sl Proof.} 
Let  $G=\langle f_{1},\ldots ,f_{s}\rangle .$ 
By \cite{sumi2}, $A(G)\cup P(G)\subset F(G)$ and 
for each $z\in \oc \setminus (A(G)\cup P(G))$, 
$\delta (f)=T_{f}(z).$ 
Let $z\in A(G)\cup P(G).$ If there exists an $n\in \Bbb{N}$ such that 
for each $\om \in \{ 1,\ldots ,s\} ^{n}$ and each $y\in f_{\om }^{-1}(z)$, 
we have $y\in \oc \setminus (A(G)\cup P(G))$, 
then by the previous argument, $T_{f}(z)=\delta (f).$ 
If there exists a strictly increasing sequence $(n_{j})_{j=1}^{\infty }$ 
of positive integers such that for each $j\in \Bbb{N}$, 
there exists an $\om \in \{ 1,\ldots ,s\} ^{n_{j}}$ and a 
$y\in f_{\om }^{-1}(z)$ with $y\in A(G)\cup P(G)$, then by 
\cite[Lemma 1.30]{sumi1}, $T_{f}(z)=\infty .$ 
Thus, we have $\delta (f)=t_{0}(f).$ We are done. 
\qed 
\end{lem}
\begin{rem}
\label{dimjexppaperrem}
Let $f=(f_{1},\ldots ,f_{s})\in \Exp(s)$ and 
let $G=\langle f_{1},\ldots ,f_{s}\rangle .$ 
Then, by \cite{sumi2} and Lemma~\ref{deltat0lem}, we have 
$\HD (J(G))\leq s_{0}(G)\leq S_{G}(z)\leq \delta (f)=T_{f}(z)=t_{0}(f),$ 
for each $z\in \oc \setminus (A(G)\cup P(G)).$ 
In addition to the above assumption, if 
$G$ satisfies the open set condition i.e., 
if there exists a non-empty open set $U$ in $\oc $ such that 
$\cup _{j=1}^{s}f_{j}^{-1}(U)\subset U$ and 
$f_{i}^{-1}(U)\cap f_{j}^{-1}(U)=\emptyset $ for each 
$(i,j)$ with $i\neq j$, then by \cite{sumi2}, 
$$\HD(J(G))=s_{0}(G)=S_{G}(z)=\delta (f)=T_{f}(z)=t_{0}(f),$$ 
for each $z\in \oc \setminus (A(G)\cup P(G)).$ 
\end{rem}
\begin{rem}
The Bowen parameter $\delta (f)$ can be strictly larger than 
two (See \cite[Example 4.14]{sumi2}). In fact, a small neighborhood $U$ of 
$(z^{2}, z^{2}/4, z^{2}/3)\in \Exp(3)$ satisfies that 
for each $f\in U$, $\delta (f)>2.$ For details, 
see Section~\ref{Examples}. 
\end{rem}
In the sequel \cite{su}, 
we will give some estimates of $\delta (f).$ 
\section{Expandingness}
In this section, we show that for an element 
$f=(f_{1},\ldots ,f_{s})\in \Exp(s)$, 
the skew product map $\tilde{f}:J(\tilde{f})\rightarrow J(\tilde{f})$ 
associated with $\{f_{1},\ldots ,f_{s}\} $ is an expanding map 
in the sense of Chapter 3 of \cite{pubook}. We need more notation.
 
 Let $f=(f_{1},\ldots ,f_{s})\in (\mbox{Rat})^{s}$ and 
 let $G:=\langle f_{1},\ldots ,f_{s}\rangle .$  
For every 
$n\le |\om|$ let $\om|_n=(\om_1,\om_2,\ld ,\om_n)$. If $\om\in\Sg_{s}^{*}$, we put
$$
[\om]=\{\tau\in\Sg_s:\tau|_{|\om|}=\om\}.
$$
If $\om,\tau\in\Sg_s\cup\Sg_s^*$, $\om\wedge\tau$ is the longest initial
subword common for both $\om$ and $\tau$. 
Let $\alpha $ be a fixed number with $0<\alpha <1/2.$ 
We endow the shift space $\Sg_s$
with the metric $\rho_\alpha $ defined as
$$
\rho_{\alpha }(\om,\tau)=\alpha ^{|\om\wedge\tau|}
$$
with the standard convention that $\alpha ^{\infty}=0$. The metric 
$d_\alpha $
induces the product topology on $\Sg_s$. Denote the spherical distance
on $\oc$ by $\hat\rho$ and equip the product
space $\Sg_s\times\oc$ with the metric $\rho$ defined as follows.
$$
\rho((\om,x),(\tau,y))=\max\{\rho_{\alpha }(\om,\tau),\hat\rho(x,y)\}.
$$
Of course $\rho$ induces the product topology on $\Sg_s\times\oc$. Using the
fiberwise expanding property (\ref{1112505}), 
\cite[Proposition 3.2]{hiroki1} and the expanding property of 
the shift map $\sg:\Sg_s\to\Sg_s$, it is fairly easy to prove the following.

\

\bthm\label{t1112505p61}
Let $G=\langle f_{1},\ldots ,f_{s}\rangle $ be an expanding 
rational semigroup. 
Let $\tilde{f}: \Sg _{s}\times \oc  \rightarrow \Sg_{s}\times \oc $ be 
the skew product associated with 
$f=(f_{1},\ldots ,f_{s}) .$ Then, 
the dynamical system $\tilde{f}:J(\tilde{f})\to J(\tilde{f})$ 
is a topologically exact open 
distance expanding map (in the sense of Chapter~3 of \cite{pubook}), 
meaning that
\begin{itemize}
\item[(a)] The map $\tilde{f}:J(\tilde{f})\to J(\tilde{f})$ is open.
\item[(b)] The map $\tilde{f}: 
J(\tilde{f})\to J(\tilde{f})$ is Lipschitz continuous.
\item[(c)] There exists $q\ge 1$ and $\d>0$ such that
$$
\rho (\tilde{f}^{q}(y),\tilde{f}^{q}(x))\ge 4\rho (y,x)
$$
for all $x,y\in J(\tilde{f})$ with 
$\rho(x,y)\le 2\d.$ 
\item[(d)] The map $\tilde{f}:J(\tilde{f})\to J(\tilde{f})$ 
is topologically exact.
\end{itemize}
In addition $q$ and $\d$ depend only on $C,\eta $ in (\ref{1112505}) 
and the Lipschitz
constant of $\tilde{f}:J(\tilde{f})\to J(\tilde{f})$. 
Note that $2\d$ is an expansive constant 
of the map $\tilde{f}:J(\tilde{f})\to J(\tilde{f})$; 
in particular it is an expansive constant
for the shift map $\sg:\Sg_s\to\Sg_s$.
\ethm

\section{$J$-stability}\label{stability}
In this section, we construct a conjugacy map $h: J(\tilde{f})
\rightarrow J(\tilde{g})$, when 
$f\in \Exp(s)$ and $g$ is close enough to $f.$ This 
conjugacy will be used to construct an analytic family of 
Perron-Frobenius operators.
 
 Define the metric $\rho _{\infty }$ on $(\mbox{Rat})^{s}$ as follows: 
for any $f=(f_{1},\ldots ,f_{s}), g=(g_{1},\ldots ,g_{s})\in (\mbox{Rat})^{s}$, we set 
$$\rho_{\infty }(f,g):=\rho _{\infty }(\tilde{f},\tilde{g})
:=\sup \{ \rho (\tilde{f}(z),\tilde{g}(z)): z\in \Sg_{s}\times \oc \} , $$ 
where $\tilde{f}$ (resp. $\tilde{g}$) denotes the 
skew product map associated with 
$f=(f_{1},\ldots ,f_{s}) $ (resp. $g=(g_{1},\ldots ,g_{s})) .$  
Given a set $D\sbt\Sg_s\times\oc$ and $r>0$, we put 
$$
B(D,r)=\{z\in\Sg_s\times\oc:\rho(z,D)<r\}.
$$
where 
$$
\rho(A,B)=\inf\{\rho(a,b):(a,b)\in A\times B\}.
$$
Similarly, given a set 
$D\subset (\mbox{Rat})^{s}$ and $r>0$, we 
put 
$$B(D,r):= \{ f\in (\mbox{Rat})^{s}\mid 
\rho _{\infty }(f,D)<r\} .$$ 

Let $\Comp^*(\Sg_s\times\oc)$ be the set of all non-empty compact (=closed) 
subsets of $\Sg_s\times\oc$ and let $\Comp^*(\oc)$ be the set of all 
non-empty compact (=closed) subsets of $\oc$. The Hausdorff metric on 
$\rho_H$ on $\Comp^*(\Sg_s\times\oc)$ is defined as follows.
$$
\rho_H(A,B)=\inf_{r>0}\{A\sbt B(B,r) \  \&  \  B\sbt B(A,r)\}.
$$
The Hausdorff metric $\hat\rho_H$ on $\Comp^*(\oc)$ is defined analogously.
Let $\Psi : \Exp(s) \rightarrow \Comp^*(\oc )$ be the 
map defined by 
$\Psi (f_{1},\ldots ,f_{s}):= J(\langle f_{1},\ldots ,f_{s}\rangle ).$ 
Then, we have the following. 
\begin{lem}
\label{jcontilem}
The map $\Psi :\Exp(s) \rightarrow \mbox{{\em Comp}}^*(\oc )$ is continuous.
\end{lem}
\sl{Proof.}
Let $f=(f_{1},\ldots ,f_{s})\in \Exp(s).$ 
By \cite[Theorem 3.1]{HM}, \cite[Lemma 2.3 (g)]{hiroki1}, 
and \cite[Lemma 3.2]{sumi2},  
we have 
\begin{equation} 
\label{jcontilemeq1} 
J(\langle f_{1},\ldots ,f_{s}\rangle )=
\overline{\{ z\in \oc : \exists u\in \langle f_{1},\ldots ,f_{s}\rangle  
\mbox{ such that }
 u(z)=z, |u'(z)|>1\} }.
\end{equation} 
Hence, for any $\epsilon >0$, there exists a finite set 
$Q_{f}=\{ \xi _{1,f},\ldots ,\xi _{l,f}\} $ such that 
$J(\langle f_{1},\ldots ,f_{s}\rangle )\subset B(Q_{f},\epsilon /2)$ and such that 
each $\xi _{j,f} $ is a repelling fixed point of some 
$u_{j,f}\in \langle f_{1},\ldots ,f_{s}\rangle .$ 
By Implicit Function Theorem, 
it follows that there exists an open neighborhood $U$ of $f$ such that 
for each $g=(g_{1},\ldots ,g_{s})\in U$ and each $j$ with $1\leq j\leq l$, 
there exists a repelling fixed point $\xi _{j,g}$ of 
some $u_{j,g}\in \langle g_{1},\ldots ,g_{s}\rangle $ such that 
$\hat{\rho }(\xi _{j,f},\xi _{j,g})\leq \epsilon /2.$ 
Therefore, setting $Q_{g}:=\{ \xi _{1,g}\ldots ,\xi _{l,g}\} $, 
we obtain that 
for each $g\in U$, 
$J(\langle f_{1},\ldots ,f_{s}\rangle )
\subset B(Q_{g},\epsilon )\subset B(J(\langle g_{1},\ldots ,g_{s}\rangle ),\epsilon
).$  
Combining this with (\ref{expopenlemeq3}) in 
the proof of Lemma~\ref{expopenlem}, 
it follows that the map $\Psi :\Exp(s)\rightarrow 
\mbox{{\em Comp}}^{\ast }(\oc )$ is continuous at 
$f.$ We are done.

\qed 
\begin{rem}
In \cite{sumihyp1}, some results which are similar to 
Lemma~\ref{expopenlem} and Lemma~\ref{jcontilem} have been 
shown, regarding the dynamics of finitely generated hyperbolic 
rational semigroups having elements of degree greater than or equal to 
two. 
\end{rem}

\bprop\label{p2112505}
The function $f=(f_{1},\ldots ,f_{s})\mapsto J(\tilde{f})$ 
from $\Exp(s)$ to $\mbox{{\em Comp}}^*(\Sg_s\times\oc)$
is continuous. 
\eprop
{\sl Proof.} Fix $f=(f_{1},\ldots ,f_{s})\in\Exp(s)$. 
Now fix $\e>0$. Since $J(\tilde{f})$ is the
closure of all repelling periodic points of $\tilde{f}$, 
there exists a finite 
set $P=\{ \xi _{1},\ldots ,\xi _{l}\} \sbt J(\tilde{f})$ of repelling periodic
points of 
$\tilde{f}$ such that
\begin{equation}
\label{3112505}
J(\tilde{f})\sbt B(P,\e/2).
\end{equation} 
Choose $p\ge 1$ so large that $\tilde{f}^{p}$ fixes every point in $P$. 
We may assume that for each $j=1,\ldots ,l$, 
$\| (\tilde{f}^{p})'(\xi _{j})\| >8.$ Since
all points $\xi_1,\xi_2,\ld,\xi_l$ of $P$ are repelling, there exists 
$0<r<\e/2$ and for every $1\le j\le l$ there exists a continuous inverse 
branch $\tilde{f}_{\xi_j}^{-p}:
B(\xi_j,4r)\to\Sg_s\times\oc$ of $\tilde{f}^p$ such that
$\tilde{f}_{\xi_j}^{-p}(\xi_j)=\xi_j$ and 
\begin{equation}
\label{2112505}
\lt| \lt|\lt(\tilde{f}_{\xi_j}^{-p}\rt)'(z)\rt|\rt| \le 1/4
\end{equation}
for all $z\in B(\xi_j,4r)$, 
where $(\tilde{f}_{\xi_j}^{-p})'(z):= 
((\tilde{f}^{p})'(\tilde{f}_{\xi _{j}}^{-p}(z)))^{-1}.$ 
 Write $\xi_j=(\ov\om^j,z_j)$, where $z_j\in\oc$
and $\ov\om^j\in\Sg_s$ is the infinite concatenation of a finite word $\om^j$
of length $p$. There then exists a unique injective meromorphic inverse
branch $f_{\om^j,j}^{-1}:B(z_j,4r)\to\oc$ of $f_{\om^j}:\oc\to\oc$ such 
that
$$
\tilde{f}_{\xi_j}^{-p}(\tau,y)=\(\om^j\tau, f_{\om^j,j}^{-1}(y)\)
$$
for all $y\in B(z_j,4r)$. In particular 
$f_{\om^j,j}^{-1}(z_j)=z_j$
and, by (\ref{2112505}), 
$$
\lt| \lt|\lt(f_{\om^j,j}^{-1}\rt)'(y)\rt|\rt|\le 1/4
$$
for all $y\in B(z_j,4r)$. Thus, there exists $\eta_1>0$ such that if 
$g=(g_{1},\ldots g_{s})\in  B(f,\eta_1)= 
\{ g=(g_{1},\ldots ,g_{s})\in \Exp(s): 
\rho _{\infty }(f,g)<\eta _{1}\} $, 
then for every $1\le j\le l$ there exists a meromorphic
inverse branch $g_{\om^j,j}^{-1}:B(z_j,2r)\to\oc$ of 
$g_{\om^j}:\oc\to\oc$ such 
that $g_{\om^j,j}^{-1}(\ov B(z_j,r))\sbt\ov B(z_j,r)$ and, furthermore,
$\|(g_{\om^j,j}^{-1})'(y)\|\le 1/2$ for all $y\in\ov B(z_j,r)$. It therefore follows
from the Banach Contraction Principle that there exists $x_j\in\ov B(z_j,r)$,
a unique fixed point of $g_{\om^j,j}^{-1}:\ov B(z_j,r)\rightarrow 
\ov B(z_j,r)$, and
$\|(g_{\om^j,j}^{-1})'(x_j)\|\le 1/2$. 
Consequently, $\tilde{g}^p(\ov\om^j,x_j)=
(\ov\om^j,x_j)$ and $\|(\tilde{g}^p)'(\ov\om^j,x_j)\| \ge 2$. Hence, 
$(\ov\om^j,x_j)\in
J(\tilde{g})$. 
Since also $\rho\((\ov\om^j,x_j),\xi_j\)=\hat\rho(x_j,z_j)<\e/2$,
using (\ref{3112505}), we get that
$$
J(\tilde{f})\sbt B(J(\tilde{g}),\e).
$$
In order to prove the "opposite" inclusion (with appropriately smaller $\eta_1$)
suppose on the contrary that there exists a sequence 
$(g_n)_{n=1}^\infty=
((g_{n,1},\ldots ,g_{n,s}))_{n=1}^{\infty }\sbt
\Exp(s)$ such that $\lim_{n\to\infty}g_n=f$ and 
$$
J(\tilde{g}_{n})\cap 
\left((\Sg_s\times\oc)\sms B(J(\tilde{f}),\e)\right)\ne\es
$$
for all $n\ge 1$. For every $n\ge 1$ choose a point $z_n$ belonging to this 
intersection. Since the space $\Sg_s\times\oc$ is compact, passing to a subsequence,
we may assume without loss of generality that the sequence 
$(z_n)_{n=1}^\infty$
converges. Denote its limit by $z$. Since 
$z\not\in B(J(\tilde{f}),\e)$, this point is 
in $F(\tilde{f})$. 
So, by Lemma~3.13, p.401 in \cite{sumi2}, there exists $q\ge 1$ 
such that $\pi_2(\tilde{f}^{q}(z))\in F(\langle f_{1},\ldots 
f_{s}\rangle )$. Applying Lemma~\ref{jcontilem}, 
we therefore conclude 
there exists $\th>0$ such that for all $n\ge 1$ large enough
$$
\pi_{2}(\tilde{f}^{q}(z))\in B(\pi_{2}(\tilde{f}^{q}(z)),\th)\sbt 
F(S_n),
$$ 
where $S_{n}:= \langle g_{n,1},\ldots ,g_{n,s}\rangle .$ 
Since $\lim_{n\to\infty}z_n=z$, we thus have that 
$$
\pi_{2}(\tilde{g}_{n}^{q}(z_n))\in
B(\pi_{2}(\tilde{f}^{q}(z)),\th)\sbt F(S_n)
$$ 
for all $n\ge 1$ large enough. On the other
hand, $\pi_{2}(\tilde{g}_{n}^{q}(z_n))\in 
\pi_{2}(J(\tilde{g}_{n}))=J(S_n)$. This contradiction 
finishes the proof. 
\qed 

\ 

Fix now $f=(f_{1},\ldots ,f_{s})\in\Exp(s)$. Then there exists $p\ge 1$
such that 
\begin{equation}
\label{1112605p158}
\|(\tilde{f}^{p})'(z)\|\ge 4
\end{equation}
for all $z\in J(\tilde{f})$. 
Since $J(\tilde{f})$ is compact and the function $z\mapsto
\|(\tilde{f}^{p})'(z)\|$ is continuous on $\Sigma _{s}\times \oc $, there exists
$\th'>0$
such that
$$ 
\|(\tilde{f}^{p})'(z)\|\ge 3
$$
for all $z\in B(J(\tilde{f}),\th')$. 
Combining this and Proposition~\ref{p2112505},
we see that that there exists $\th''\in(0,\th']$ such that
$$ 
\|(\tilde{g}^{p})'(z)\|\ge 2
$$
for all $g=(g_{1},\ldots ,g_{s})\in B(f,\th'')$ and all 
$z\in B(J(\tilde{g}),\th')$.

 Now using the above, in particular Proposition~\ref{p2112505}, fairly
straightforward continuity type considerations and \cite[Theorem 2.14]{sumi1} lead to the following.

\

\blem\label{l2112605p158}
Suppose that $f=(f_{1},\ldots ,f_{s})\in\Exp(s)$ 
and let $p\ge 1$ be given by (\ref{1112605p158}).
Then there exists a number $\th=\th_{f}^{(1)}>0$ 
such that the following properties 
are satisfied.
\begin{itemize}
\item[(a)] For all $g=(g_{1},\ldots ,g_{s})\in 
B(f,\th_{f}^{(1)})$, all $x\in\ov B(J(\tilde{g}),\th_{f}^{(1)})$
and all $y\in\ov B(x,\th_{f}^{(1)})$,
$$
\|(\tilde{g}^{p})'(y)\|\ge 2, \  \rho(\tilde{g}^{p}(x),\tilde{g}^{p}(y))
\ge 2\rho(x,y)
$$
and $\tilde{g}^{p}|_{B(x,\th_{f}^{(1)})}$ 
is one-to-one and $g^{p}\(B(x,\th_{f}^{(1)})\)
\spt B(\tilde{g}^{p}(x),2\th_{f}^{(1)})$. 
\item[(b)] If $g\in B(f,\th_{f}^{(1)})$ and $\{\tilde{g}^{pn}(x):n\ge 0\}\sbt 
B(J(\tilde{g}),\th_f^{(1)})$, then $x\in J(\tilde{g})$.
\end{itemize}
\elem

\

\ni As a direct consequence of item (a) of this lemma we get the following.

\

\blem\label{l3112605p159}
Suppose that $f=(f_{1},\ldots ,f_{s})\in\Exp(s)$ 
and let $p\ge 1$ be given by (\ref{1112605p158}).
If $g=(g_{1},\ldots ,g_{s})\in B(f,\th_{f}^{(1)})$, $n\ge 1$, and 
$\tilde{g}^{pk}(x)\in\ov B(J(\tilde{g}),\th_{f}^{(1)})$
for all $0\le k\le n$, then there exists a unique continuous inverse branch
$\tilde{g}_{x}^{-pn}:
B\(\tilde{g}^{pn}(x),\th_{f}^{(1)}\)\to 
B(x,\th_{f}^{(1)})$ of $\tilde{g}^{pn}$
sending $\tilde{g}^{pn}(x)$ to $x$. 
In addition, $\tilde{g}_{x}^{-pn}$ is Lipschitz continuous
with Lipschitz constant $\le 2^{-n}$ and 
$\tilde{g}_{x}^{-pn}\(B\(\tilde{g}^{pn}(x),
\th_{f}^{(1)}\)\)\sbt B(x,\th_{f}^{(1)})$.
\elem

\ni From now onwards, unless otherwise stated, assume that the integer $p\ge 1$ 
ascribed to $f$ by (\ref{1112605p158}) is equal to $1$. We then call $f$ simple.
\ 

 Recall that a sequence $(x_i)_{i=0}^n\sbt \Sg_s\times\oc$, $0\le n\le\infty$ is 
called a $\g$-pseudoorbit with respect to the map 
$\tilde{f}:\Sg_s\times\oc\to
\Sg_s\times\oc$ provided that
$$
\rho(\tilde{f}(x_i),x_{i+1})\le\g
$$
for all $0\le i\le n-1$. The pseudoorbit 
$(x_i)_{i=0}^n$ is said to be $\b$-shadowed
by a point $x\in \Sg_s\times\oc$ provided that
$$
\rho(x_i,\tilde{f}^i(x))\le\b
$$
for all $0\le i\le n$.

\ni We shall prove the following. 

\

\blem\label{l4112605p159}
Assume the same as in Lemma~\ref{l2112605p158} (and so the same as in 
Lemma~\ref{l3112605p159}). Fix $\b\in(0,\th_{f}^{(1)}/2]$. Suppose that 
$(x_i)_{i=0}^n\sbt B(J(\tilde{g}),\th_{f}^{(1)}/4)$ is a 
$\b$-pseudoorbit for the 
skew product map $\tilde{g}$. 
For every $0\le i\le n-1$, let $y_i=\tilde{g}_{x_{i}}^{-1}(x_{i+1})$,
which is defined since $x_i,\tilde{g}(x_i)\in B(J(\tilde{g}),\th_{f}^{(1)})$ 
and 
$x_{i+1}\in B\(\tilde{g}(x_i),\b\)\sbt B\(\tilde{g}(x_i),\th_{f}^{(1)})$. 
Then for all
$0\le i\le n-1$, we have
\begin{itemize}
\item[(a)] $y_i\in B(x_i,\th_{f}^{(1)}/2)$ and $\tilde{g}(y_i)=x_{i+1}\in 
B(J(\tilde{g}),\th_{f}^{(1)})$. So, 
in view of Lemma~\ref{l3112605p159}, each 
inverse branch $\tilde{g}_{y_i}^{-1}:B\(\tilde{g}(y_i),\th_{f}^{(1)}\)\to 
B(y_i,\th_{f}^{(1)}/2)$ is well defined.
\item[(b)] For all $0\le i\le n-1$
$$
\tilde{g}_{y_i}^{-1}\(\ov B(x_{i+1},\b)\)\sbt \ov B(x_i,\b)
$$
and, consequently, all the compositions 
$$
\overline{g}_{i}^{-i}:=\tilde{g}_{y_0}^{-1}\circ 
\tilde{g}_{y_1}^{-1}\circ \cdots \circ 
\tilde{g}_{y_{i-1}}^{-1}:
\ov B(x_i,\b)\to\Sg_s\times\oc
$$
are well defined for all $i=1,2,\ld,n$.
\item[(c)] $\(\overline{g}_{i}^{-i}(\ov B(x_i,\b))\)_{i=0}^n$ 
is a descending sequence 
of non-empty compact sets.
\item[(d)] $\bi_{i=0}^{n}\overline{g}_{i}^{-i}
(\ov B(x_i,\b))\ne\es$ and all the 
elements of this intersection $\b$-shadow the pseudoorbit $(x_i)_{i=0}^n$.
\item[(e)] If $n=+\infty$, then the intersection in the item (d) is a
singleton which belongs to $J(\tilde{g})$.
\end{itemize}
\elem
{\sl Proof.} Since $(x_i)_{i=0}^n$ is a $\b$-pseudoorbit, 
$\tilde{g}(x_i)\in
B(x_{i+1},\b)\sbt B(J(\tilde{g}),\th_{f}^{(1)}/2)\sbt 
B(J(\tilde{g}),\th_{f}^{(1)})$. By the 
Lipschitz
part of Lemma~\ref{l3112605p159}, $y_i=\tilde{g}_{x_i}^{-1}(x_{i+1})\in
B(x_i,\b/2)\sbt B(x_i,\th_{f}^{(1)}/2)\sbt B(J(\tilde{g}),\th_{f}^{(1)})$ 
and item 
(a) is proved.
In order to prove item (b) take an arbitrary point $z\in 
\ov B(x_{i+1},\b)$, $0\le i\le n-1$. Applying the Lipschitz
part of Lemma~\ref{l3112605p159} again, we get
$$
\aligned
\rho\(\tilde{g}_{y_i}^{-1}(z),x_i\)
&\le \rho\(\tilde{g}_{y_i}^{-1}(z),y_i\)+\rho(y_i,x_i) \\
&=   \rho\(\tilde{g}_{y_i}^{-1}(z),\tilde{g}_{y_i}^{-1}(x_{i+1})\)+
     \rho\(\tilde{g}_{x_i}^{-1}(x_{i+1}),\tilde{g}_{x_i}^{-1}
     (\tilde{g}(x_i))\)\\
&\le 2^{-1}\rho(z,x_{i+1})+2^{-1}\rho(x_{i+1},\tilde{g}(x_i))\\
&\le 2^{-1}\b+2^{-1}\b
 =\b.
\endaligned
$$
Hence, $\tilde{g}_{y_i}^{-1}\(\ov B(x_{i+1},\b)\)\sbt \ov B(x_i,\b)$
and item (b) is proved. Item (c) is now an immediate consequence of
(b), and the first part of (d) is an immediate consequence of
(c) and compactness of the sets $\overline{g}_{i}^{-i}(\ov B(x_i,\b))$. Since
for every $0\le k\le n$,
$$
\tilde{g}^{k}\lt(\bi_{i=0}^n\overline{g}_{i}^{-i}
(\ov B(x_i,\b))\rt)
\sbt \tilde{g}^k(\overline{g}_{k}^{-k})(\ov B(x_k,\b))
=\ov B(x_k,\b),
$$
the second part of item (d) follows. Since $\bi_{i=0}^n
\overline{g}_{i}^{-i}(\ov B(x_i,\b))=\overline{g}_{n}^{-n}
(\ov B(x_n,\b))$, it follows
from Lemma~\ref{l3112605p159} that 
$\diam\lt(\bi_{i=0}^{n}\overline{g}_{i}^{-i}(\ov B(x_i,\b))\rt)\le 2^{-n}\b$,
and the singleton part of item (e) follows. Obviously
\begin{equation}\lab{1112905p160}
\overline{g}_{n}^{-n}(\ov B(x_n,\b))\sbt \overline{g}_{n}^{-n}
(\ov B(x_n,\th_{f}^{(1)}/2))
\end{equation}
and $\ov B(x_n,\th _{f}^{(1)}/2)\cap J(\tilde{g})\ne\es$ as $x_n\in 
B(J(\tilde{g}),\th_{f}^{(1)}/4)$.
Since the set $J(\tilde{g})$ is completely invariant under $\tilde{g}$, we
conclude that $J(\tilde{g})\cap \overline{g}_{n}^{-n}
(\ov B(x_n,\th_{f}^{(1)}/2))\ne\es$.
Thus $\bi_{n=0}^\infty \overline{g}_{n}^{-n}
(\ov B(x_n,\th_{f}^{(1)}/2))$ is a 
singleton belonging to $J(\tilde{g})$. The second part of item (e) is then
concluded by invoking (\ref{1112905p160}). 

\qed 

\

\ni As a straightforward consequence of Lemma~\ref{l2112605p158} 
and Lemma~\ref{l3112605p159} we get the following.

\

\bprop\label{p1112905p160}
Assume the same as in Lemma~\ref{l2112605p158}. Then for every 
$g\in B(f,\th_f^{(1)})$, the number $\th_f^{(1)}$ is an expansive
constant of $\tilde{g}:J(\tilde{g})\rightarrow J(\tilde{g})$ 
meaning that if $x,y\in J(\tilde{g})$ and 
$\rho(\tilde{g}^n(y),\tilde{g}^n(x))\le \th_f^{(1)}$ for all 
$n\ge 0$, then $x=y$.
\eprop 

\

\ni Now if $y$ and $z$ $\b$-shadow the same pseudoorbit $(x_n)_{n=0}^\infty$,
then $\rho(\tilde{g}^n(y),\tilde{g}^n(z))\le 2\b\le\th_f^{(1)}$. Thus, as an 
immediate consequence of Lemma~\ref{l4112605p159} and 
Proposition~\ref{p1112905p160}, we get the first part of the following.

\

\bprop\label{p2112905p161} (shadowing lemma)
Assume the same as in Lemma~\ref{l2112605p158}. Then for every 
$g\in B(f,\th_f^{(1)})$ and every $\b\in (0,\th_f^{(1)}/2]$,
every $\b$-pseudoorbit $(x_n)_{n=0}^\infty\sbt B(J(\tilde{g}),\th_f^{(1)}/4)$ 
is $\b$-shadowed by a unique element $x\in J(\tilde{g})$. If in addition
$\sigma ^{n}(\om )=\pi_1(x_n)$ for all $n\ge 0$, then $\pi_1(x)=\om$.
\eprop

\ni In order to see the second part of this proposition, just notice 
that $\rho_\alpha (\sg^n(\pi_1(x)),\sg^n(\om))\le\b\le \th_f^{(1)}$ and 
$\th_f^{(1)}$ is an expansive constant for the shift map
$\sg:\Sg_s\to\Sg_s$. 

By Proposition~\ref{p2112505} there exist
$\th_f^{(2)}\in (0,\th_f^{(1)}/3]$ such that 
\begin{equation}\lab{3112905p161}
\rho_H(J(\tilde{g}),J(\tilde{f}))<\th_f^{(1)}/4
\end{equation}
whenever $\rho _{\infty }(g,f)\le \th_f^{(2)}$. We shall prove the following 
main result of this section.

\

\bthm\label{t3112905p161}
Suppose that $f\in\Exp(s)$ is simple. Then for every $g\in 
B(f,\th_f^{(2)})$ there exists a unique homeomorphism $h=h_g:
J(\tilde{f})\to J(\tilde{g})$ with the following properties.
\begin{itemize}
\item[(a)] $\tilde{g}\circ h=h\circ \tilde{f}$,
\item[(b)] $\pi_1\circ h=\pi_1$,
\item[(c)] The homeomorphism $h:J(\tilde{f})\to J(\tilde{g})$ is H\"older
continuous with the H\"older exponent $\ka_f$ and the same 
H\"older constant $L_f$, which is thus the same for all $g\in 
B(f,\th_f^{(2)})$.
\item[(d)] $\rho_{\infty, J(\tilde{f})}(h,\Id)
:= \sup \{ \rho (h(z),Id(z)): z\in J(\tilde{f})\} 
\le \th_f^{(1)}/2$.
\end{itemize}
\ethm
{\sl Proof.} Fix $g\in B(f,\th_f^{(2)})$. Fix also $z\in J(\tilde{f})$. 
Then using (\ref{3112905p161}), we get $(\tilde{f}^n(z))_{n=0}^\infty\sbt
B(J(\tilde{g}),\th_f^{(1)}/4)$ and 
$$
\rho\(\tilde{f}^{n+1}(z),\tilde{g}(\tilde{f}^n(z))\)
=\rho\(\tilde{f}(\tilde{f}^n(z)),\tilde{g}(\tilde{f}^n(z))\)
\le\rho_\infty(\tilde{f},\tilde{g})
<\th_f^{(2)}
\le\th_f^{(1)}/3.
$$
Therefore, in view of Proposition~\ref{p2112905p161} (shadowing lemma)
applied with $\b=\th_f^{(1)}/3$, there exists a unique element, denote
it by $h(z)\in J(\tilde{g})$, that $\th_f^{(1)}/3$-shadows the pseudoorbit
$(\tilde{f}^n(z))_{n=0}^\infty$. 
In addition $\pi_1(h(z))=\pi_1(z)$. So, we have
defined a map $h:J(\tilde{f})\to J(\tilde{g})$ 
such that in particular, the property
(b) is satisfied. Also, for every $n\ge 0$,
\begin{equation}\lab{4112905p162}
\rho\(\tilde{g}^n(h(z)),\tilde{f}^n(z)\)\le\b=\th_f^{(1)}/3,
\end{equation}
and reading it with $n=0$, we get item (d). Also, reading (\ref{4112905p162})
for all $n\ge 1$, in the form $\rho\(\tilde{g}^{n-1}(\tilde{g}(h(z))),
\tilde{f}^{n-1}(\tilde{f}(z))\)\le \th_f^{(1)}/3$, 
we see that the point $\tilde{g}(h(z))$ \ 
$\th_f^{(1)}/2$-shadows the pseudoorbit 
$(\tilde{f}^n(\tilde{f}(z)))_{n=0}^\infty$. Thus
$\tilde{g}(h(z))=h(\tilde{f}(z))$ and item (a) is established. In order to prove 
item (c) fix $L\ge 1$, a Lipschitz constant of $\tilde{f}:\Sg_s\times\oc\to
\Sg_s\times\oc$. Take $x,y\in J(\tilde{f})$ with $0<\rho(x,y)<\th_f^{(1)}/3$.
Let $k\ge 0$ be the largest integer such that
\begin{equation}\lab{5112905p162}
L^k\rho(x,y)<\th_f^{(1)}/3.
\end{equation}
Then
\begin{equation}\lab{6112905p162}
L^{k+1}\rho(x,y)\geq \th_f^{(1)}/3.
\end{equation}
It follows from (\ref{5112905p162}) that $\rho(\tilde{f}^j(x),\tilde{f}^j(y))<
\th_f^{(1)}/3$ for all $j=0,1,2,\ld,k$. Hence, invoking (\ref{4112905p162}),
we get for every $j=0,1,2,\ld,k$ that
$$
\aligned
\rho\(\tilde{g}^j(h(x)),\tilde{g}^j(h(y))\)
&\le \rho\(\tilde{g}^j(h(x)),\tilde{f}^j(x)\)
+\rho\(\tilde{f}^j(x),\tilde{f}^j(y)\)+
    \rho\(\tilde{f}^j(y),\tilde{g}^j(h(y))\) \\
&<{1\over 3}\th_f^{(1)}+{1\over 3}\th_f^{(1)}+{1\over 3}\th_f^{(1)} \\
&=\th_f^{(1)}.
\endaligned
$$
Hence, Lemma~\ref{l3112605p159} and Lemma~\ref{l2112605p158} yield
$$
\aligned
\rho(h(x),h(y))
&=\rho\(\tilde{g}_{h(x)}^{-k}(\tilde{g}^k(h(x))),
\tilde{g}_{h(x)}^{-k}(\tilde{g}^k(h(y)))\) \\
&\le 2^{-k}\rho\(\tilde{g}^k(h(x))),\tilde{g}^k(h(y))\) \\
&\le 2^{-k}\th_f^{(1)}.
\endaligned
$$
Since, by (\ref{6112905p162}), $2^{-k}=L^{-k{\log 2\over\log L}}
\le (3L/\th_f^{(1)})^{{\log 2\over\log L}}
\rho^{{\log 2\over\log L}}(x,y)$, we therefore conclude that
$$
\rho(h(x),h(y))
\le \th_f^{(1)}(3L/\th_f^{(1)})^{{\log 2\over\log L}}
\rho^{{\log 2\over\log L}}(x,y),
$$
and the H\"{o}lder continuity in condition (c) is proved. In order to prove that
$h:J(\tilde{f})\to J(\tilde{g})$ is $1$-to-$1$ 
suppose that $x,y\in J(\tilde{f})$ and
$h(x)=h(y)$. Then we get from (\ref{4112905p162}) that 
$\rho(\tilde{f}^n(x),\tilde{f}^n(y))\le 2\th_f^{(1)}/3$ for all $n\ge 0$, equality
$x=y$ follows from Proposition~\ref{p1112905p160}. So, in order to
complete the proof of our theorem, we only need to show that 
$h:J(\tilde{f})\to J(\tilde{g})$ is surjective. Indeed, fix $y\in J(\tilde{g})$.
Reasoning analogously as in the very beginning of the proof, we get that
$$
(\tilde{g}^n(y))_{n=0}^\infty\sbt B(J(\tilde{f}),\th_f^{(1)}/4) \  \text{ and }
\rho\(\tilde{g}^{n+1}(y),\tilde{f}(\tilde{g}^n(y))\)\le\th_f^{(1)}/3
$$
for all $n\ge 0$. Therefore, again analogously as in the beginning of 
the proof, we conclude that there exists a point $x\in J(\tilde{f})$ that
$\th_f^{(1)}/3$-shadows the pseudoorbit $(\tilde{g}^n(y))_{n=0}^\infty$. 
Writing out what this means, we get that $\rho(\tilde{f}^n(x),\tilde{g}^n(y))<
\th_f^{(1)}/3$. But this also means that the 
$(\th_f^{(1)}/3)$-pseudoorbit $(\tilde{f}^n(x))_{n=0}^\infty$ is 
$\th_f^{(1)}/3$-shadowed (with respect to the map $\tilde{g}$ 
by the point $y$).
So, $y=h(x)$ and we are done. 
\qed 
\

\section{Perron-Frobenius Operators and Their Regularity
Properties}
\label{Perron}
In this section, we investigate some complex analytic families of 
Perron-Frobenius operators on the Banach space of H\"older continuous 
functions.. 

 Fix a simple element $f\in\Exp(s)$. 
Denote by $C_0$ the 
Banach space of all complex-valued continuous functions on $J(\tilde{f})$
endowed with the supremum norm $\| \cdot  \| _{\infty }$. Fix $\a>0$. Given 
$\varphi \in C_0$ let
$$
V_\a=\inf\{L\ge 0:|\varphi (y)-\varphi (x)|\le
L\rho(y,x)^\a \  \text{ for all } \ x,y\in J(\tilde{f})\},
$$
be the $\a$-variation of the function $\varphi $ and let
$$
||\varphi ||_\a=V_\a(\varphi )+||\varphi ||_\infty.
$$
Clearly the space
$$
\H_\a=\H_\a(J(\tilde{f}))=\{\varphi \in C_{0} :||\varphi ||_\a<\infty\}
$$
endowed with the norm $||\cdot||_\a$ is a Banach space densely
contained in $C_0$ with respect to the $||\cdot||_\infty$ norm. 
Each member of $\H_\a$ is called a H\"older continuous function
with exponent $\a$. 
For each Banach space $B$, we denote by $L(B)$ the space of bounded 
operators on $B$ endowed with the operator norm. 

Let $W$ be an open subset of $\C^q$ with 
some $q\ge 1$. Suppose that for every $\l\in W$ there is given a 
function $\zeta_\l\in\H_\a$. The formula
$$
\L_\l \varphi (z)=\sum_{x\in \tilde{f}^{-1}(z)}e^{\zeta_\l(x)}\varphi (x)
$$
defines a bounded linear operator acting on the Banach space $C_0$, 
called the Perron-Frobenius operator of the potential $\zeta_\l$. It 
is straightforward to check that $\L_\l$ preserves the Banach 
space $\H_\a$. We shall prove the following. 

\

\blem\lab{l1013106p167} 
 If the map $\l\mapsto \L_{\l}\in L(\H_\a)$, 
$\l\in W$, is continuous and for every $z\in J(\tilde{f})$ the function
$\l\mapsto \zeta_\l(z)$, $\l\in W$, is holomorphic, then the map
$\l\mapsto\L_{\l}\in L(\H_\a)$ is holomorphic on $W$.
\elem
{\sl Proof.} It follows directly from our assumptions that for every
$\varphi \in H_\a$ and every $z\in J(\tilde{f})$ the function
$$
\l\mapsto \L_{\l}\varphi (z)=
\sum_{x\in \tilde{f}^{-1}(z)}e^{\zeta _{\lambda }(x)}\varphi (x)\in 
{\C}, \  \l\in W,
$$
is holomorphic. Fix a one-dimensional round disk $D\sbt W$.
Let $\gamma\subset D$ be a simple closed curve. Fix $\varphi \in H_\a$ 
and $z\in J(\tilde{f})$. 
By Cauchy's theorem $\int_\gamma \L_{\l}\varphi (z)d\l=0$.
Since the function $\l\mapsto \L_{\l}\varphi \in H_\a$ 
is continuous, the integral
$\int_\gamma\L_{\l} \varphi  d\l$ exists and for every $z\in J(\tilde{f})$, 
we have
$(\int_\gamma\L_{\l}\varphi  d\l)(z)=\int_\gamma\L_{\l}\varphi (z)d\l=0$.
Hence, $\int_\gamma\L_{\l}\varphi  d\l=0$. Now, since the function 
$\l\mapsto
\L_{\l}\in L(H_\a)$ is continuous, the integral
$\int_\gamma\L_{\l}d\l$ exists and for every $\varphi \in H_\a$,
$(\int_\gamma\L_{\l}d\l)(\varphi )=\int_\gamma\L_{\l}\varphi d\l=0$.
Thus, $\int_\gamma\L_{\l}d\l =0$ and in view of Morera's theorem, 
the function $\l\mapsto\L_{\l}\in L(H_\a)$ is holomorphic in $D$. 
So, we are done in virtue of Hartogs Theorem. The proof is complete.
\qed 

\

\noindent {\bf Notation:} 
For any parameter $\lambda _{0}\in \Bbb{C}^{q}$ and 
any $R>0$, we set 
$B(\lambda _{0},R):= \{ \lambda \in \Bbb{C}^{q}: \| \lambda -\lambda _{0}\| 
<R\} .$ 

\ 

\ni In order for Lemma~\ref{l1013106p167} to be applicable we shall prove the 
following.

\

\blem\lab{l2013106p167} 
If for every $\l\in W$, we have that $\zeta_\l\in\H_\a$, 
the function $\l\mapsto \zeta_\l\in\H_\a$ is continuous and for every 
$z\in J(\tilde{f})$ the function $\l\mapsto \zeta_\l(z)$, $\l\in W$, 
is holomorphic, then the map $\l\mapsto\L_{\l}\in L(\H_\a)$, 
$\l\in W$, is continuous.
\elem

{\sl Proof.} Fix $\l_0\in W$ and $R>0$ so small that $B(\l_0,R)
\sbt W$. By our assumptions there exists $M>0$ such that
\begin{equation}\lab{2013106p168}
|\zeta_\l(z)|\le M
\end{equation}
for all $z\in J(\tilde{f})$ and all $\l\in B(\l_0,R/2)$.
Let 
$$
\hat M:=\sup\lt\{{|e^t-1|\over |t|}:t\in \ov B(0,2M)\rt\}<\infty.
$$
It suffices to show that the function $\l\mapsto\L_{\l}\in 
L(\H_\a)$ is Lipschitz continuous with the same Lipschitz constant
separately with respect to each variable. Therefore, taking a
complex plane parallel to the coordinate planes,
we may assume without loss of generality that $q=1$. Put 
$D=B(\l_0,R/4)$. It follows from Cauchy's formula that 
\begin{equation}\lab{3013106p168}
|\dot\zeta_\l(z)|\le 16MR^{-1}
\end{equation}
for all $z\in J(\tilde{f})$ and all 
$\l\in D$, where $\dot\zeta_\l(z)
={d\over d\l}\zeta_\l(z)$. Hence, for all $\l,\l'\in D$ and 
$x\in J(\tilde{f})$ we obtain
$$
|\zeta_{\l'}(x)-\zeta_\l(x)|
=\lt|\int_\l^{\l'}\dot\zeta_\mu(x)d\mu \rt|
\le 16MR^{-1}|\l'-\l|.
$$
Therefore, 
\begin{equation}\lab{3020206p168}
\aligned
|e^{\zeta_{\l'}(x)}-e^{\zeta_\l(x)}|
&=   \exp(\re(\zeta_\l(x)))\lt|\exp\(\zeta_{\l'}(x)-\zeta_\l(x)\)-1\rt|
 \le \hat Me^M|\zeta_{\l'}(x)-\zeta_\l(x)| \\
&\le 16M\hat Me^MR^{-1}|\l'-\l|.
\endaligned
\end{equation}
Consequently, for all $\varphi \in \H_\a$,
$$
\aligned
|(\L_{\l'}-\L_\l)\varphi (z)|
&=   |\L_{\l'}\varphi (z)-\L_\l \varphi (z)|
 =   \lt|\sum_{x\in \tilde{f}^{-1}(z)}
  \varphi (x)\(e^{\zeta_{\l'}(x)}-e^{\zeta_\l(x)}\)\rt|\\
&\le ||\varphi ||_\infty
\sum_{x\in \tilde{f}^{-1}(z)}|e^{\zeta_{\l'}(x)}-e^{\zeta_\l(x)}|\\
&\le 16M\hat Me^MR^{-1}\deg(\tilde{f})||\varphi ||_\infty|\l'-\l|.
\endaligned
$$
Hence, 
\begin{equation}\lab{4020206p168}
||(\L_{\l'}-\L_\l)\varphi ||_\infty
\le 16M\hat Me^MR^{-1}\deg(\tilde{f})|\l'-\l|||\varphi ||_\a
\end{equation}
for all $\l,\l'\in D$ and  all $\varphi \in \H_\a$. In order to estimate
the $\a$-variation of the function $(\L_{\l'}-\L_\l)\varphi $, 
let $\lambda '\in D$ be a point and consider the
function $\psi:D\times J(\tilde{f})\times J(\tilde{f})\to\C$ 
given by the formula
$$
\psi(\l,w,z)
=e^{\zeta_\l(w)}-e^{\zeta_\l(z)}+e^{\zeta_{\l'}(z)}-e^{\zeta_{\l'}(w)}.
$$
Obviously
\begin{equation}\lab{2020206p168}
\psi(\l',w,z)=0.
\end{equation}
Put $\dot\psi(\l,w,z)={\bd\over\bd\l}\psi(\l,w,z)$. Applying
(\ref{2013106p168}) and (\ref{3013106p168}), we get
\begin{equation}\lab{4013106p169}
\aligned
|\dot\psi(\l,w,z)|
&=   \lt|{\bd\over\bd\l}e^{\zeta_\l(w)}-{\bd\over\bd\l}e^{\zeta_\l(z)}\rt|
 =   \lt|e^{\zeta_\l(w)}\dot\zeta_\l(w)-e^{\zeta_\l(z)}\dot\zeta_\l(z)\rt| \\
&=   \lt|e^{\zeta_\l(w)}\(\dot\zeta_\l(w)-\dot\zeta_\l(z)\)
        +\dot\zeta_\l(z)\(e^{\zeta_\l(w)}-e^{\zeta_\l(z)}\)\rt| \\
&\le e^M|\dot\zeta_\l(w)-\dot\zeta_\l(z)|
        +16MR^{-1}|e^{\zeta_\l(w)}-e^{\zeta_\l(z)}|.
\endaligned
\end{equation}
Since the function $\l\mapsto \zeta_\l\in\H_\a$ is continuous, we 
have that $V=\sup\{V_\a(\zeta_\l):\l\in B(\lambda _{0},R/2)\}<+\infty$. 
Applying Cauchy's formula, we have for all $z\in J(\tilde{f})$ 
and all $\l\in D$
that
$$
\aligned
|\dot\zeta_\l(w)-\dot\zeta_\l(z)|
&=   \lt|{1\over 2\pi i}\int_{\bd B(\l _{0}, R/2)}
{\zeta_\xi(w)-\zeta_\xi(z)\over (\xi-\l)^2}d\xi\rt|
 \le {16\over 2\pi R^2}\int_{\bd B(\l _{0},R/2)}|\zeta_\xi(w)-\zeta_\xi(z)||d\xi| \\
&\le 16R^{-1}V_\a(\zeta_\xi)\rho (w,z)^\a
 \le 16R^{-1}V\rho (w,z)^\a.
\endaligned
$$
We can therefore continue as follows.
\begin{equation}\lab{1020206p169}
\aligned
|\dot\psi(\l,w,z)|
&\le 16R^{-1}Ve^M\rho (w,z)^\a+16MR^{-1}e^M\hat M|\zeta_\l(w)-\zeta_\l(z)| \\
&\le \(16R^{-1}Ve^M+16MR^{-1}e^M\hat MV)\rho (w,z)^\a.
\endaligned
\end{equation}
Put $T=16VR^{-1}e^M(1+M\hat M)$. Applying (\ref{1020206p169}) and using 
(\ref{2020206p168}), we get for all \nl $(\l,w,z)\in 
D\times J(\tilde{f})\times 
J(\tilde{f})$ that
$$
\aligned
|\psi(\l,w,z)|
   &=|\psi(\l,w,z)-\psi(\l',w,z)|
    =\lt|\int_{\l'}^\l\dot\psi(\mu,w,z)d\mu\rt| \\
 &\le\int_{\l'}^\l|\dot\psi(\mu,w,z)||d\mu|
\le T|\l'-\l|\rho (w,z)^\a.
\endaligned
$$
Therefore, using also (\ref{3020206p168}), we obtain for all $x,y\in
J(\tilde{f})$ with $\rho (x,y)<\theta _{f}^{(1)}$ that
$$
\aligned
|(\L_\l &-\L_{\l'})\varphi (y)-(\L_\l-\L_{\l'})\varphi (x)|=\\
&=  \lt|\sum_{z\in \tilde{f}^{-1}(y)}
 \(e^{\zeta_\l(z)}-e^{\zeta_{\l'}(z)}\)\varphi (z)
       -\sum_{w\in \tilde{f}^{-1}(x)}
       \(e^{\zeta_\l(w)}-e^{\zeta_{\l'}(w)}\)\varphi (w)\rt|\\
&=\bigg|  \sum 
[\(e^{\zeta_\l(z)}-e^{\zeta_{\l'}(z)}\)(\varphi (z)-\varphi (w)) \\
       &  \  \  \  \  \  \  \   \  \ -\(e^{\zeta_\l(w)}-e^{\zeta_\l(z)}
              +e^{\zeta_{\l'}(z)}-e^{\zeta_{\l'}(w)}\)\varphi (w)]\bigg|\\
&=\lt|  \sum 
  \(e^{\zeta_\l(z)}-e^{\zeta_{\l'}(z)}\)(\varphi (z)-\varphi (w))
       -\psi(\l,w,z)\varphi (w)\rt|\\
&\le   \sum  
    |e^{\zeta_\l(z)}-e^{\zeta_{\l'}(z)}|\cdot|\varphi (z)-\varphi (w)|
       +|\psi(\l,w,z)|\cdot|\varphi (w)|\\
&\le    16M\hat Me^{M}R^{-1}\deg(\tilde{f})|\l-\l'|||\varphi ||_\a\rho(y,x)^\a 
        +T||\varphi ||_\infty\deg(\tilde{f})|\l-\l'|\rho(y,x)^\a \\
&\le    (16M\hat Me^{M}R^{-1}+T)\deg(\tilde{f})|\l-\l'|||\varphi ||_\a
\rho(y,x)^\a .
\endaligned
$$
Thus, combining this with (\ref{4020206p168}), 
$$
V_\a\((\L_\l-\L_{\l'})\varphi \)
\le (32M\hat Me^{M}R^{-1}+T+(\theta _{f}^{(1)})^{\alpha })\deg(\tilde{f})|\l-\l'|||\varphi ||_\a,
$$
and combining this with (\ref{4020206p168}) again,  
we get $$||(\L_\l-\L_{\l'})\varphi ||_\a
\le \deg(\tilde{f})|(48M\hat{M}e^{M}R^{-1}+T+(\theta _{f}^{(1)})^{\alpha })|\l-\l'|||\varphi ||_\a.$$ So,
$$
||\L_\l-\L_{\l'}||\le \deg(\tilde{f})|
(48M\hat{M}e^{M}R^{-1}+T+(\theta _{f}^{(1)})^{\alpha })|\l-\l'|.
$$
We are done. 

\qed 

\

\section{Analytic Extensions of Perron-Frobenius Operators.}
\label{Analytic}
In this section, for a given analytic family of 
expanding rational semigroups, we construct an associated 
real-analytic family of Perron-Frobenius operators, and then 
we provide a complex analytic extension of this family, 
in order to use the results proven in the previous section.   

 Let us first describe in detail the setting of this section. $\La$ is 
assumed to be an open subset of a finite dimensional complex Banach space (ex.
$\C^d$). 

\sp\ni Let $\{G_\l=\langle f_{\l ,1},\ldots ,f_{\l ,s}\rangle \}_{\l\in\La}$ 
be an analytic family of
expanding rational semigroups. 
For every $\l\in\La$ put $f_\l=(f_{\l ,1},\ldots ,f_{\l ,s})
\in (\mbox{{\em Rat}})^{s}$. 
Fix $\l_0\in\La$ and put $f=f_{\l _{0}}$.

Then we easily obtain the following. 
\

\bprop\lab{p4112905p163}
The map $\l\mapsto f_\l$ is continuous.
\eprop

\

For each $\l \in \Lambda $, let 
$\theta _{f_{\lambda }}^{(i)}$ $(i=1,2)$ be 
the number for $f_{\lambda }$ obtained in Section~\ref{stability}. 
\ni So, for every $\l\in\La$ there exists $R_\l>0$ so small that $f_\g\in 
B\(f_\l,\frac{1}{2}\min\{\th_{f_{\l}}^{(1)},\th_{f_{\l}}^{(2)}\}\)$ whenever $\g\in
B(\l,R_\l)$. 
For each $\l \in \Lambda $, let 
$\tilde{f}_{\l }:\Sigma _{s}\times \oc \rightarrow 
\Sigma _{s}\times \oc $ be the skew product map associated with 
$f_{\l }=(f_{\l, 1},\ldots ,f_{\l ,s}) .$  
We shall prove the following.

\

\blem\lab{l5112905p164}
If $\l_0\in\La$ and for every $\l\in B(\l_0,R_{\l_0})$, 
$h_\l:J(\tilde{f})\to
J(\tilde{f}_\l)$ is the unique conjugating homeomorphism coming from 
Theorem~\ref{t3112905p161}, then for every $z\in J(\tilde{f})$ the map 
$\l\mapsto\pi_2(h_\l(z))\in\oc$, $\l\in B(\l_0,R_{\l_0})$, is holomorphic.
\elem
{\sl Proof.} 
Conjugating $\langle f_{\l _{0},1},\ldots ,f_{\l _{0}, s}\rangle $ by a M\"{o}bius 
map, we may assume that $\pi _{2}(J(\tilde{f}))=
J(\langle f_{\lambda _{0},1},\ldots ,f_{\lambda _{0},s}\rangle )
\subset \Bbb{C} .$ By Lemma~\ref{jcontilem}, 
we may assume that 
$J(\langle f_{\l ,1},\ldots ,f_{\l ,s}\rangle )\subset \Bbb{C}$ 
for each $\l \in \Lambda $, 
in order to prove our lemma. 

Set $R_0=R_{\l_0}$. Fix a repelling periodic point 
$(\om,x)\in J(\tilde{f})$ of
$\tilde{f}$, say of period $p\ge 1$. Consider the map $H(z,\l)=
\pi_2(\tilde{f}_\l^p(\om,z))-z$, 
$(z,\l)\in\C\times\La$. Then $H(x,\l_0)=0$ and
$$
{\bd\over\bd z}H(z,\l)|_{(x,\l_0)}=(\tilde{f}^p)'(\om,x)-1\ne 0
$$
since $|(\tilde{f}^p)'(\om,x)|>1$. It therefore follows from the Implicit 
Function Theorem that there exists $\hat R(\om,x)>0$ and a holomorphic 
map $u_{\om,x}:B(\l_0,\hat R(\om,x))\to\C$ such that $u_{\om,x}(\l_0)
=x$ and $H(u_{\om,x}(\l),\l)=0$ for all $\l\in B(\l_0,\hat R(\om,x))$.
But $\pi_2(h_{\l_0}(\om,x))=\pi_2(\om,x)$ and, because of 
Theorem~\ref{t3112905p161},
$$
\aligned
H(\pi_2(h_\l(\om,x)),\l)
&=\pi_2\circ \tilde{f}_\l^p(\om,\pi_2(h_\l(\om,x)))-\pi_2(h_\l(\om,x))\\
&=\pi_2\circ \tilde{f}_\l^p(h_\l(\om,x))-\pi_2(h_\l(\om,x)) \\
&=\pi_2\circ h_\l \circ \tilde{f}^p(\om,x)-\pi_2(h_\l(\om,x)) \\
&=\pi_2(h_\l(\om,x))-\pi_2(h_\l(\om,x))
=0.
\endaligned
$$
Therefore, in view of the uniqueness part of the Implicit 
Function Theorem, there exists $R(\om,x)\in(0,\hat R(\om,x)]$ such
that $u_{\om,x}(\l)=\pi_2(h_\l(\om,x))$ for all $\l\in 
B(\l_0,R(\om,x))$. In particular, the map $\l\mapsto
\pi_2(h_\l(\om,x))$, $\l\in B(\l_0,R(\om,x))$ is holomorphic. Now suppose
that the map $\l\mapsto \pi_2(h_\l(\om,x))$ defined on 
$B(\l_0,R_0)$ fails to be holomorphic. Select then a parameter
$\l_1\in B(\l_0,R_0)$ such that the map $\l\mapsto 
\pi_2(h_\l(\om,x))$ fails to be holomorphic at $\l_1$ but it is 
holomorphic throughout $B(\l_0,||\l_1-\l_0||)$. Obviously, $||\l_1-\l_0||
\ge R(\om,x)>0$. Since $h_{\l_1}(\om,x)=\lim_{\l\to\l_1}h_\l(\om,x)$ 
(the limit taken throughout $B(\l_0,||\l_1-\l_0||)$), we have that
$\tilde{f}_{\l_1}^p\(h_{\l_1}(\om,x)\)=h_{\l_1}(\om,x)$; also 
$|(\tilde{f}_{\l _{1}}^p)'(h_{\l_1}(\om,x))|>1$ since $f_{\l_1}\in\Exp(s)$. 
Replacing in the above considerations $\l_0$ by $\l_1$ and $x$ by
$\pi_2(h_{\l_1}(\om,x))$ thus yields that the map $\l\mapsto
\pi_2\(h_\l^1(h_{\l_1}(\om,x))\)$ is holomorphic on a neighborhood
of $\l_1$, where $h_\l^1:J(\tilde{f}_{\l_1})\to J(\tilde{f}_\l)$ 
is the unique 
conjugating homeomorphism coming from Theorem~\ref{t3112905p161}.
Since $h_\l(\om,x)=h_\l^1\circ h_{\l_1}(\om,x)$ on this neighborhood,
we thus conclude that the map $\l\mapsto\pi_2(h_\l(\om,x))$ 
is holomorphic on a neighborhood of $\l_1$, which is a contradiction.
We have thus proved the following.

\sp\ni {\bf Claim~1.} For every periodic point $\xi\in J(\tilde{f})$
the map $\l\mapsto\pi_2(h_\l(\xi))\in\oc$, $\l\in 
B(\l_0,R_0)$, is holomorphic.

\sp\ni Now fix an arbitrary point $z_\infty\in J(\tilde{f})$ and 
let $(z_n)_{n=1}^\infty$ be a sequence or repelling periodic points 
of $\tilde{f}$ converging to $z_\infty$. Define the maps $u_n:
B(\l_0,R_0)\to\oc$, $n=1,2,\ld,\infty$, by the formula
$u_n(\l)=\pi_2(h_\l(z_n))$. By Claim~1 all these maps with finite
$n$ are holomorphic. It follows from Theorem~\ref{t3112905p161}(c)
that there exists a constant $L _{\lambda _{0}}$ and a 
constant $\kappa _{\l _{0}}$ such that 
for all $1\le n<\infty$ and all $\l\in B(\l_0,R_0)$, 
$$\hat{\rho }(u_n(\l),u_\infty(\l))
=\hat{\rho }(\pi_2(h_\l(z_n)),\pi_2(h_\l(z_\infty)))
\le\rho(h_\l(z_n),h_\l(z_\infty))
\le L_{\l_0}\rho(z_n,z_\infty)^{\ka_{\l_0}}.
$$
Hence, the sequence $(u_n)_{n=1}^\infty$ of holomorphic maps 
converges uniformly on $B(\l_0,R_0)$ to the map
$u_\infty$. So, $u_\infty:B(\l_0,R_0)\to\oc$ is holomorphic, 
and the proof is complete. 
\qed 

\begin{rem}
Theorem~\ref{t3112905p161} and Lemma~\ref{l5112905p164} imply that 
for each $\om \in \Sigma _{s}$, 
$J_{\om }(\tilde{f}_{\lambda })$ moves by a 
holomorphic motion (\cite{Mc}).
\end{rem}

\

\ni Now suppose that our Banach space containing $\La$ is equal 
to $\C^d$ with some $d\ge 1$. Embed $\C^d$ into $\C^{2d}$ by 
the formula
$$
(x_1+iy_1,x_2+iy_2,\ld,x_d+iy_d)\mapsto (x_1,y_1,x_2,y_2,\ld,x_d,y_d).
$$
For every $z\in\C^d$ and every $r>0$ denote by $D_d(z,r)$ the 
$d$-dimensional polydisk in $\C^d$ centered at $z$ and with "radius" $r$.
We will need the following lemma, which is of general dynamics 
independent character.

\

\blem\lab{l1120105p165} 
For every $M\ge 0$, for every $R>0$, for every $\l^0\in\C^d$, and for 
every analytic function $\phi:D_d(\l^0,R)\to\C$ bounded in modulus by 
$M$ there exists an analytic function $\^\phi:D_{2d}(\l^0,R/4)\to\C$
that is bounded in modulus by $4^dM$ and whose restriction to the
polydisk $D_d(\l^0,R/4)$ coincides with $\re\phi$, the real part of 
$\phi$.
\elem
{\sl Proof.} Denote by $\N_0$ the set of all non-negative integers. 
Write the analytic function $\phi:D_d(\l^0,R)\to\C$ in the form of 
its Taylor series expansion
$$
\phi(\l_1,\l_2,\ld,\l_d)
=\sum_{\a\in\N_0^d}a_\a(\l_1-\l_1^0)^{\a_1}(\l_2-\l_2^0)^{\a_2}\ld
                       (\l_d-\l_d^0)^{\a_d}.
$$
By Cauchy's estimates we have
\begin{equation}\lab{1013106p166}
|a_\a|\le {M\over R^{|\a|}}
\end{equation}
for all $\a\in\N_0^d$. We have
$$
\aligned
\re\phi(\l_1 &,\l_2,\ld,\l_d)=\\
&=\sum_{\a\in\N_0^d}\re\bigg[a_\a\lt(\sum_{p=0}^{\a_1}
      {\a_1\choose p}\(\re\l_1-\re\l_1^0\)^p
                     \(\im\l_1-\im\l_1^0\)^{\a_1-p}i^{\a_1-p}\rt)\cdot \\
& \  \  \  \  \  \  \  \  \cdot\lt(\sum_{p=0}^{\a_2}
      {\a_2\choose p}\(\re\l_2-\re\l_2^0\)^p 
                    \(\im\l_2-\im\l_2^0\)^{\a_2-p}i^{\a_2-p}\rt)\cdot\ld \\
& \  \  \  \  \  \  \  \  \ld\cdot\lt(\sum_{p=0}^{\a_d}
      {\a_1\choose p}\(\re\l_d-\re\l_d^0\)^p
                     \(\im\l_d-\im\l_d^0\)^{\a_d-p}i^{\a_d-p}\rt)\bigg]\\
&=\sum_{\b\in\N_0^{2d}}\re\lt[a_{\hat\b}\prod_{j=1}^d
      {\b_j^{(1)}+\b_j^{(2)}\choose\b_j^{(1)}}i^{\b_j^{(2)}}
      \(\re\l_j-\re\l_j^0\)^{\b_j^{(1)}}\(\im\l_j-\im\l_j^0\)^{\b_j^{(2)}}\rt]\\
&=\sum_{\b\in\N_0^{2d}}\re\lt(a_{\hat\b}\prod_{j=1}^d
      {\b_j^{(1)}+\b_j^{(2)}\choose\b_j^{(1)}}i^{\b_j^{(2)}}\rt)
      \(\re\l_j-\re\l_j^0\)^{\b_j^{(1)}}\(\im\l_j-\im\l_j^0\)^{\b_j^{(2)}},
\endaligned
$$
where we wrote $\b\in\N_0^{2d}$ in the form 
$\(\b_1^{(1)},\b_1^{(2)},\b_2^{(1)},\b_2^{(2)},\ld,\b_d^{(1)},\b_d^{(2)}\)$ 
and we also put $\hat\b=
\(\b_1^{(1)}+\b_1^{(2)},\b_2^{(1)}+\b_2^{(2)},\ld,\b_d^{(1)}+\b_d^{(2)}\)\in\N_0^d$.
Set
$$
c_\b=\re\lt(a_{\hat\b}\prod_{j=1}^d
      {\b_j^{(1)}+\b_j^{(2)}\choose\b_j^{(1)}}i^{\b_j^{(2)}}\rt).
$$
Using (\ref{1013106p166}), we get
$$
|c_\b|
\le |a_{\hat\b}|\prod_{j=1}^d{\b_j^{(1)}+\b_j^{(2)}\choose\b_j^{(1)}}
\le MR^{-|\hat\b|}\prod_{j=1}^d2^{\b_j^{(1)}+\b_j^{(2)}}
=   MR^{-|\b|}2^{|\b|}.
$$
Thus the formula
$$
\^\phi(x_1,y_1,x_2,y_2,\ld,x_d,y_d)
=\sum_{\b\in\N_0^{2d}}c_\b\prod_{j=1}^d
 \(x_j-\re\l_j^0\)^{\b_j^{(1)}}\(y_j-\im\l_j^0\)^{\b_j^{(2)}}
$$
defines an analytic function on $D_{2d}(\l_0,R/4)$ and
$$
|\^\phi(x_1,y_1,x_2,y_2,\ld,x_d,y_d)|\le 4^dM.
$$
Obviously $\^\phi|_{D_d(\l_0,R/4)}=\re\phi|_{D_d(\l_0,R/4)}$, and 
we are done. 

\qed 

\

\ni Let $\{ G_{\l }=\langle f_{\l ,1},\ldots ,f_{\l,s}\rangle \} 
_{\l \in \Lambda }$ be an analytic family of expanding rational semigroups.  
Coming back to the dynamical situation (and keeping the assumption that
$\La\sbt\C ^{d}$), 
we assume that 
for each $\l \in \Lambda $, 
$J(G_{\lambda })\subset \Bbb{C}.$ 
\begin{rem}
\label{assumerem}
Note that in order to prove the main results of this paper, 
we may assume the above. For, 
by Lemma~\ref{jcontilem},  
for a fixed $\l _{0}\in \Lambda $,  
there exists a neighborhood $B$ of $\l _{0}$ and a 
M\"{o}bius transformation $u $ 
such that 
for each $\l \in B$, $J(u G_{\l }u ^{-1})\subset \Bbb{C}$, 
where $u G_{\l }u ^{-1}:=  
\{ u gu ^{-1}: g\in G_{\l }\} .$ 
\end{rem}
Fix an element $\l _{0}\in \Lambda $ and set 
$f=f_{\lambda _{0}}$, set 
$f_{\l }:=(f_{\l ,1},\ldots ,f_{\l ,s})$, etc., and 
we use the notation which was given in the beginning of this section. 
For every $z\in J(\tilde{f})$ consider the function
\begin{equation}\lab{3020306p12}
\l\mapsto\psi_z(\l):={\tilde{f}_\l'(h_\l(z))\over \tilde{f}'\(z)}, 
\ \l \in D_{d}(\l _{0},R_{0}),\ R_{0}:= R_{\l _{0}}.
\end{equation}
It follows from Lemma~\ref{l5112905p164} 
that this function is holomorphic. 
Taking $R_0>0$ and $\th_{f_{\l_0}}^{(1)}>0$
sufficiently small, we have that
$$
\lt|{\tilde{f}_\l'(w)\over \tilde{f}'\(z)}-1\rt|<1/6
$$
whenever $\rho(w,z)<\th_{f_{\l_0}}^{(1)}$ and $\l \in D_{d}(\l _{0},R_{0})$. 
It 
therefore follows from Theorem~\ref{t3112905p161}(d) that
\begin{equation}\lab{2020306p12}
|\psi_z(\l)-1|<1/5
\end{equation}
for all $z\in J(\tilde{f})$ and all $\l\in D_{d}(\l_0,R_0)$. 
Hence, taking a branch $\log : \{ z\in \C: |z-1|<1/5\} \rightarrow \C $ 
of logarithm with $\log (1)=0$, for every $z\in J(\tilde{f})$
there exists $\log\psi_z:D_d(\l_0,R_0)\to\C$, a unique holomorphic branch of
logarithm of $\psi_z$ sending $\l_0$ to $0$. It is bounded in modulus by $1/4$.
Now, in virtue of Lemma~\ref{l1120105p165}, with $R_*\in (0,R_0)$, there 
exists $\^{\re\log\psi_z}:D_{2d}(\l_0,R_*)\to\C$, an analytic extension of 
$\re\log\psi_z:D_d(\l_0,R_0)\to\R$ bounded in modulus by $4^{d-1}$. Now for 
all $(t,\l,z)\in \C\times D_{2d}(\l_0,R_*)\times J(\tilde{f})$, put 
\begin{equation}\lab{8020306p13}
\zeta_{(t,\l)}(z)=-t\^{\re\log\psi_z}(\l)+t\log|\tilde{f}'(z)|
\end{equation}
Of course for all $z\in J(\tilde{f})$ all the maps 
$(t,\l)\mapsto\zeta_{(t,\l)}(z)$
are holomorphic. Let 
$$
\ka:=\ka_{f_{\l_0}}
$$
coming from Theorem~\ref{t3112905p161}(c). Aiming to apply 
Lemma~\ref{l1013106p167}, we shall prove the following.

\

\blem\lab{l1020206p13}
There exists an $\tilde{R}\ (<R_{\ast })$ such that 
for every $(t,\l)\in \C\times D_{2d}(\l_0,\tilde{R})$, the function
$\zeta_{(t,\l)}:J(\tilde{f})\to\C$ is in $\H_\ka $ and the map $(t,\l)\mapsto
\zeta_{(t,\l)}\in \H_\ka$ is continuous.
\elem
{\sl Proof.} Note that for all $(t,\l)\in \C\times D_{2d}(\l_0,R_*)$, 
\begin{equation}\lab{5020206p13}
||\zeta_{(t,\l)}||_\infty\le 4^{d-1}|t|+|t|\cdot || \log|\tilde{f}'|||_\infty.
\end{equation}
Fix $r_1>0$ so small that all the maps $\tilde{f}_{\l}$, $\l\in 
D_d(\l_0,R_0)$, $1\le i\le s$, are univalent on
all disks centered at points of $J(\tilde{f}_\l)$ with radii equal to $r_{1}$.
Starting with (\ref{2020306p12}), applying Lemma~\ref{l5112905p164} and 
making use of item (c) of Theorem~\ref{t3112905p161}, we deduce that 
there exists a constant $M_{1}\geq 1$ and a constant $0<R_{1}<R_{\ast }$ such that 
\begin{equation}
\label{psirhoeq}
|\psi _{z}(\l )-\psi _{w}(\l )|\leq M_{1}\rho (z,w)^{\kappa }
\end{equation}
for all $w,z\in J(\tilde{f})$ and all $\l \in D_{d}(\l _{0},R_{1}).$ 
Hence, there exists a constant $M_{2}\geq 1$ such that 
\begin{equation}
\label{logpsirhoeq}
|\log \psi _{z}(\l )-\log \psi _{w}(\l )|\leq M_{2}\rho (z,w)^{\kappa }
\end{equation}
for all $w,z\in J(\tilde{f})$ and all $\l \in D_{d}(\l _{0},R_{1}).$ 
Applying Lemma~\ref{l1120105p165}, we get 
\begin{equation}
\label{relogpsirhoeq}
|\^{\re\log \psi _{z}}(\l )-\^{\re\log \psi _{w}}(\l )|\leq 4^{d}M_{2}
\rho (z,w)^{\kappa }
\end{equation}
for all $w,z\in J(\tilde{f})$ and all $\l \in D_{2d}(\l _{0},R_{1}/4).$ 
Thus, there exists a constant $M>0$ such that
$$
|\zeta_{(t,\l)}(w)-\zeta_{(t,\l)}(z)|\le M|t|\rho (w,z)^\ka
$$
for all $w,z\in J(\tilde{f})$ and all 
$(t,\l)\in \C\times D_{2d}(\l_0,R_1/4)$. Along
with (\ref{5020206p13}) this shows that 
$\zeta_{(t,\l)}\in\H_\ka$ for all $(t,\l )\in 
\C \times D_{2d}(\l _{0},R_{1}/4).$ Now we pass
to the continuity part of the proof. 
Since $\^{\re\log\psi_z}:D_{2d}(\l_0,R_{\ast })\to\C$ is bounded
in modulus by $4^{d-1}$, taking $R_{1}$ small enough, by Cauchy's formula 
we conclude that for all $z\in J(\tilde{f})$ 
and all $\l,\l'\in
D_{2d}(\l_0,R_1/4)$, we have
$$
|\^{\re\log\psi_z}(\l')-\^{\re\log\psi_z}(\l)|\le T\|\l'-\l\|
$$
with some universal constant $T>0$. Consequently, 
for 
all $z\in J(\tilde{f})$ and all $\l,\l'\in
D_{2d}(\l_0,R_{1}/4)$ and all $t\in\C$,
\begin{equation}\lab{1020306p14}
||\zeta_{(t,\l')}-\zeta_{(t,\l)}||_\infty\le |t|T\|\l'-\l\|.
\end{equation}
Also, for all $z\in J(\tilde{f})$, 
all $\l\in D_{2d}(\l_0,R_1/4)$ and all $t_0,t\in\C$, 
we have
$$
|\zeta_{(t,\l)}(z)-\zeta_{(t_0,\l)}(z)|
\le |t-t_0||\^{\re\log\psi_z}(\l)|+|t-t_0|\| \log|f_0'|\|_\infty
\le (4^{d-1}+\| \log|f_0'|\|_\infty)|t-t_0|.
$$
As a direct consequence of this inequality and (\ref{1020306p14}), we 
conclude that
\begin{equation}\lab{7020306p14}
\aligned
\forall t_0\in\C\  \exists L_1>0\  \forall \l,\l' &\in D_{2d}(\l_0 ,R_1/4)\  
\forall t,t'\in D_{1}(t_0,2) \\
&||\zeta_{(t',\l')}-\zeta_{(t,\l)}||_\infty\le L_1||(t',\l')-(t,\l)||.
\endaligned
\end{equation}
In order to cope with the H\"older variation fix $\l\in D_{2d}(\l_0 ,R_1/4)$
and $t,t'\in\C$. Then for all $w,z\in J(\tilde{f})$ we have
$$
\aligned
|(\zeta_{(t,\l)}(w)-\zeta_{(t',\l)}(w)) 
&-(\zeta_{(t,\l)}(z)-\zeta_{(t',\l)}(z))|=\\
&=   \lt|(t-t')(\^{\re\log\psi_z}(\l)-\^{\re\log\psi_w}(\l))\rt| \\
&\le |t-t'||\^{\re\log\psi_z}(\l)-\^{\re\log\psi_w}(\l)| \\
&\le 4^dM_{2}|t-t'|\rho (w,z)^\ka.
\endaligned
$$
Hence
\begin{equation}\lab{6020306p14}
V_\ka\(\zeta_{(t,\l)}-\zeta_{(t',\l)}\)\leq 4^dM_{2}|t-t'|.
\end{equation}
In order to consider $V_{\kappa }(\zeta _{(t,\l )}-\zeta _{(t,\l ')})$, 
by (\ref{relogpsirhoeq}) and Cauchy's formula, it follows that there exists a constant $M_{3}\geq 1$ 
such that 
\begin{equation}
\label{relogpsirhoeq2}
|(\^{\re\log\psi _{z}}(\l )-\^{\re\log\psi _{w}}(\l ))
-(\^{\re\log\psi _{z}}(\l ')-\^{\re\log\psi _{w}}(\l '))|\leq 
M_{3}\rho (z,w)^{\kappa }\| \l -\l '\| 
\end{equation}
for all $w,z\in J(\tilde{f})$ and all $\l,\l' \in D_{2d}(\l _{0},R_{1}/8).$ 
Hence, we get for all $t_{0}\in \C $ and all 
$(z,w,t,\l,\l')\in J(\tilde{f})^{2}\times D_{1}(t_{0},2)\times 
D_{2d}(\l _{0},R_{1}/8)$, 
$$|(\zeta _{(t,\l )}-\zeta _{(t,\l ')})(z)-
(\zeta _{(t,\l )}-\zeta _{(t,\l ')})(w)|\leq 
(|t_{0}|+2)M_{3}\rho (z,w)^{\kappa }\| \l-\l'\|.$$ 
So $V_{\kappa }(\zeta _{(t,\l )}-\zeta _{(t,\l ')})\leq 
(|t_{0}|+2)M_{3}\| \l -\l '\| $, 
and invoking (\ref{6020306p14}), we deduce that 
for any $t_{0}\in \C $ there exists a constant $L_{2}\geq 1$ 
such that for any 
$(t, t', \l,\l ')\in D_{1}(t_{0},2)^{2}\times D_{2d}(\l _{0},R_{1}/8)^{2}$, 
$$V_{\kappa } (\zeta _{(t,\l )}-\zeta _{(t',\l ')})
\leq L_{2}\| (t,\l )-(t',\l')\| .$$ 
And bringing up (\ref{7020306p14}), we finally get 
$$\| \zeta _{(t,\l )}-\zeta _{(t',\l')}\| _{\kappa }\leq 
(L_{1}+L_{2})\| (t,\l)-(t',\l')\|$$ 
for all $(t,t',\l,\l')\in D_{1}(t_{0},2)^{2}\times 
D_{2d}(\l _{0},R_{1}/8)^{2}.$ 

We are done. 
\qed

\

\section{Real Analyticity of Bowen parameters and Hausdorff Dimensions}
\label{Real}
This section is devoted to prove our main results by applying the 
tools developed in the previous sections. 

 We now prove Theorem A. 

\ 

{\sl Proof of Theorem A:} 
 Let $\Lambda $ be a finite dimensional complex manifold. 
Let $\{G_\l\}_{\l\in\La}$, where $G_{\l }= 
\langle f_{\l ,1},\ldots ,f_{\l s}\rangle $, be an analytic family of 
expanding rational
semigroups. In order to prove Theorem A, 
as in Remark~\ref{assumerem}, 
we may assume that $\Lambda $ is an open subset of $\C ^{d}$ and 
that for each $\l \in \Lambda $, 
$J(G_{\l })\subset \C .$ 
 
Fix $\l_0\in\La$ and 
we use the notations 
$f,f_{\l}, $ etc., which were given in Section~\ref{Analytic}. 
Consider the family of
potentials $\phi_{(t,\l)}:J(\tilde{f})\to\R$, 
$(t,\l )\in\R\times D_d(\l_0,R_0)$
given by the formula
$$
\phi_{(t,\l)}(w)=-t\log|\tilde{f}_\l'(h_\l(w))|.
$$
Let $\P_\l(t)$ be the topological pressure of the potential $\phi_{(t,\l)}$
with respect the map $\tilde{f}:J(\tilde{f})\to J(\tilde{f})$, 
and let $\L_{(t,\l)}:C(J(\tilde{f}))\to
C(J(\tilde{f}))$ be the Perron-Frobenius operator defined by the formula
\begin{equation}\lab{9020306p16}
\L_{(t,\l)}\varphi (z)=\sum_{x\in \tilde{f}^{-1}(z)}
e^{\phi_{(t,\l)}(x)}\varphi (x)
               =\sum_{x\in \tilde{f}^{-1}(z)}
               |\tilde{f}_\l'(h_\l(x))|^{-t}\varphi (x).
\end{equation}
Note that $\P_\l(t)=\tilde{\P}_\l(t)$, where 
$\tilde{\P}_\l(t)$ is the
topological pressure of the potential $\^\phi_{(t,\l)}(w)=
-t\log|\tilde{f}_\l'(w)|$
with respect to the map $\tilde{f}_\l:J(\tilde{f}_\l)\to J(\tilde{f}_\l)$. 
Since we are assuming $J(G_{\l })\subset \C $ and 
the Euclidian metric and the spherical metric are comparable 
on a compact subset of $\C $, we see that 
$\tilde{\P}_{\l }(t)$ is equal to the topological 
pressure $P(t,f_{\l })$ of potential 
$-t\log \| \tilde{f}_{\l }'\| $ with respect to 
the map $\tilde{f}_{\l }:J(\tilde{f}_{\l })\rightarrow 
J(\tilde{f}_{\l })$, where 
$\| \cdot \| $ denotes the norm of the derivative with respect to 
the spherical metric on $\oc .$ 
Hence, by \cite{sumi2}, for every $\l\in D_d(\l_0,R_{0})$, 
the function $t\mapsto \P _{\l }(t)$ is continuous and 
strictly decreasing from $+\infty $ to $-\infty $, and has a unique zero 
$\delta (f_{\l })$, which is the Bowen parameter of 
$f_{\l }=(f_{\l, 1},\ldots ,f_{\l ,s}).$  

Note that $\phi_{(t,\l)}
=\zeta_{(t,\l)}$ for all $(t,\l)\in\R\times D_d(\l_0,R_{\ast })$, where 
$\zeta_{(t,\l)}$ are the potentials defined by (\ref{8020306p13}). In 
particular the formula
$$
\L_{(t,\l)}\varphi (z)=\sum_{x\in \tilde{f}^{-1}(z)}e^{\zeta_{(t,\l)}(x)}\varphi (x)
$$
extends (\ref{9020306p16}) to $\C\times D_{2d}(\l_0,R_{\ast })$. In view of
Lemma~\ref{l1020206p13}, analyticity of the maps $(t,\l)\mapsto 
\zeta_{(t,\l)}(z)$, Lemma~\ref{l2013106p167} and Lemma~\ref{l1013106p167}, the
function $(t,\l)\mapsto\L_{(t,\l)}\in L(\H_\ka)$, $(t,\l)\in\C\times 
D_{2d}(\l_0,\tilde{R})$, is holomorphic. One of the central facts (see \cite{pubook}) of the 
thermodynamic formalism of distance expanding mappings and H\"older
continuous potentials, adapted to our setting, is that $e^{\P_\l(t)}$ is a
simple isolated eigenvalue of the operator $\L_{(t,\l)}:\H_\ka\to\H_\ka$
for all $(t,\l)\in\R\times D_d(\l_0,\tilde{R})$. It therefore follows from 
the theory of perturbations for linear operators 
(\cite{Ka}, Kato-Rellich Theorem, p212, p368-370)
that given $t_0\in\R$ there exists $R\in(0,\tilde{R})$ and an analytic 
function $\g:D_1(t_0,R)\times D_{2d}(\l_0,R)\to\C$ such that for every
$(t,\l)\in D_{1}(t_{0},R)\times D_{2d}(\l_0,R)$, 
$\g(t,\l)$ is a simple isolated 
eigenvalue of the operator $\L_{(t,\l)}:\H_\ka\to\H_\ka$ and $\g(t_0,\l_0)
=e^{\P_{\l_0}(t_0)}$. Using continuity of the function $(t,\l)\mapsto
\P_\l(t)$, $(t,\l)\in (t_0-R,t_0+R)\times D_d(\l_0,R)$, we easily deduce, 
decreasing $R>0$ if necessary, that $\g(t,\l)=e^{\P_\l(t)}$ for all
$(t,\l)\in (t_0-R,t_0+R)\times D_d(\l_0,R)$. Thus 
we obtain the following: \\ 
{\bf Claim:} 
The function $(t,\l)
\mapsto\P_\l(t)$, $(t,\l)\in (t_0-R,t_0+R)\times D_d(\l_0,R)$, is 
real-analytic. 

Since 
\begin{equation}
\label{pderivative}
{\bd \P_\l(t)\over \bd t}=-\int\log|(\tilde{f}_\l)'\circ h_\l|d\mu_{(t,\l)}<0,
\end{equation}
where $\mu_{(t,\l)}$ is the Gibbs (equilibrium) state of the potential 
$\phi_{(t,\l)}$ with respect the map $\tilde{f}:J(\tilde{f})\to 
J(\tilde{f})$ (see \cite[Chapter 4, Theorem 4.6.5]{pubook}), applying the 
Implicit Function Theorem, we get that 
$\l \mapsto \delta (f_{\l })$ is real-analytic.

 We now prove the rest of Theorem A. 
 We may assume that $\Lambda $ is an open subset of $\Bbb{C}.$ 
 Moreover, as before, we may assume that 
 $J(G_{\l })\subset \Bbb{C}$ for each $\l \in \Lambda .$  
Let $\lambda _{0}\in \Lambda $ be a point.  
We now prove that $\l \mapsto 1/\delta (f_{\l})$ is 
(pluri)superharmonic around $\l _{0}.$   
If $s=1$ and $\deg (f_{\l _{0}, 1})=1$, then $\delta (f_{\l })\equiv 0$ 
around $\l _{0} .$ Hence, 
we may assume that either (1)$s>1$ or (2)$s=1$ and $\deg (f_{\l _{0},1})>1.$ 
In both cases, we  have $\delta (f_{\l _{0}})>0.$ Hence, 
we may assume that for each $\l \in \Lambda $, 
$\delta (f_{\l })>0.$ 
We set $f=f_{\l _{0}}$ and use the same notations as before. 
By the variational principle, 
we have that for each $\l $, 
\begin{equation}
\label{variationaleq1}
0=P_{\l }(\delta (f_{\l}))=\sup _{\mu \in  M(\tilde{f})}
\left\{ h_{\mu }(\tilde{f})-
\delta (f_{\l })\int _{J(\tilde{f})}\log |(\tilde{f}_{\l })'\circ h_{\l }|
d\mu \right\} , 
\end{equation}
where $ M(\tilde{f})$ denotes the set of all 
$\tilde{f}$-invariant Borel probability measures on $J(\tilde{f}).$ 
From this formula, 
we obtain 
\begin{equation}
\label{variationaleq2}
\delta (f_{\l })=\sup _{\mu \in  M(\tilde{f})}
\frac{h_{\mu }(\tilde{f})}{\int _{J(\tilde{f})}
\log |(\tilde{f}_{\l })'\circ h_{\l }|d\mu }.
\end{equation} 
So, 
\begin{equation}
\label{variationaleq3}
\frac{1}{\delta (f_{\l })}=
\inf _{\mu \in  M(\tilde{f})}
\frac{\int _{J(\tilde{f})}
\log |(\tilde{f}_{\l })'\circ h_{\l }|d\mu }
{h_{\mu }(\tilde{f})}.
\end{equation}
Therefore, the function 
$\l \mapsto 1/\delta (f_{\l })$ is an infimum of a family of 
positive (pluri)harmonic functions. 
Hence, $\l \mapsto 1/\delta (f_{\l} )$ is (pluri)superharmonic. 
Thus, we have proved that 
$\l \mapsto 1/\delta (f_{\l })$ is plurisuperharmonic. 

 Since the functions $x\mapsto 1/x$ and $x\mapsto -\log x$ are  
 decreasing convex functions in $(0,\infty )$, 
 Jensen's inequality implies that 
$\l \mapsto \delta (f_{\l })$ and $\l \mapsto \log \delta (f_{\l })$ 
are plurisubharmonic.

 Finally, let $t\in \Bbb{R}$ be a fixed number. 
 By the variational principle, we have 
 $$P_{\l }(t)=\sup _{\mu \in M(\tilde{f})}\{ h_{\mu }(\tilde{f})-
 t\int _{J(\tilde{f})}\log | (\tilde{f}_{\l })'\circ h_{\l }| \ d\mu \} .$$ 
 Hence the function 
 $\l \mapsto P_{\l }(t)$ is equal to the supremum of a  
family of pluriharmonic functions of $\l \in \Lambda .$ 
Therefore, $\l\mapsto P_{\l }(t)$ is plurisubharmonic. 
We are done. 
\qed 

\


\begin{cor}
\label{thAcor0}
Under the assumption of Theorem A, for each $\l \in \Lambda $, 
let 
$\mu _{\l }$ be the maximal entropy measure of 
$\tilde{f}_{\l }:J(\tilde{f}_{\l })\rightarrow J(\tilde{f}_{\l })$ and 
let $\tau _{\l}$ be the Gibbs (equilibrium) state of 
the potential $-\delta (f_{\l })\log \| \tilde{f}_{\l }'\| $ 
with respect to the map 
$\tilde{f}_{\l }: J(\tilde{f}_{\l })\rightarrow J(\tilde{f}_{\l }).$ 
Then, the functions 
$\l \mapsto \int _{J(\tilde{f}_{\l })}\log \| \tilde{f}_{\l} '\| 
d\mu _{\l },\ 
\l \mapsto \int _{J(\tilde{f}_{\l })}\log \| \tilde{f}_{\l }'\| 
d\tau _{\l},\ 
\l \mapsto h_{\tau _{\l}}(\tilde{f}_{\l }),$ where $ \l \in \Lambda $, 
are real-analytic.  

\end{cor}

{\sl Proof.} 
As in the proof of Theorem A, we may assume that 
for each $\l \in \Lambda $, $J(G_{\l })\subset \Bbb{C}.$ 
We use the same notation as that in the proof of Theorem A. 
As in the proof of Theorem A, we have that 
$\P_{\l }(t)=\tilde{\P}_{\l }(t)=P(t,f_{\l })$ and 
that $(t,\l )\mapsto P(t, f_{\l })$ is real-analytic. 
Since $\l \mapsto \delta (f_{\l })$ is real-analytic, 
the formula (\ref{pderivative}) and 
the equation $\delta (f_{\l })=
h_{\tau _{\l }}(\tilde{f}_{\l })/\int _{J(\tilde{f}_{\l })}
\log \| \tilde{f}_{\l }'\| d\tau _{\l }$ 
imply that the statement of 
the theorem holds. We are done.

\qed 

\

\begin{cor}
\label{thAcor1}
Under the assumption of Theorem A, suppose $\Lambda \subset \Bbb{C}.$ 
Let $\varphi (\l ):= \delta (f_{\l }).$ Then, 
$\varphi \triangle \varphi \geq 2|\nabla \varphi |^{2}$ 
in $\Lambda .$ 
\end{cor}
{\sl Proof.} 
As in the proof of Theorem A, we may assume that 
$\varphi (\l )>0$ for each $\l \in \Lambda .$ 
Since $\varphi $ is real-analytic and $1/\varphi $ is superharmonic, 
the inequality $\triangle (1/\varphi )\leq 0$ implies the 
inequality of the corollary.
\qed 

\begin{cor}
\label{thAcor2}
Under the assumption of Theorem A, 
suppose $\Lambda $ is connected and 
let 
$\varphi (\l )=\delta (f_{\l }).$ 
If $\varphi $ is pluriharmonic in a non-empty open subset  
of $\Lambda $, then $\varphi $ is constant in $\Lambda .$ 
\end{cor}
{\sl Proof.} 
If $\varphi $ is pluriharmonic in a non-empty subdomain $U$ of $\Lambda $, 
then Corollary~\ref{thAcor1} implies that 
$\varphi $ is constant in $U.$ Since $\varphi $ is real-analytic, 
it follows that $\varphi $ is constant in $\Lambda .$ 
\qed 

\

We now prove Theorem B. 

\ 

{\sl Proof of Theorem B:} 
 The first author has proved in \cite{sumi2} that if 
$G=\langle f_{1},\ldots ,f_{s}\rangle $ 
is an expanding
rational semigroup satisfying the open set condition 
i.e. there exists a non-empty open subset 
$U$ of $\oc $ such that 
$\cup _{j=1}^{s}f_{j}^{-1}(U)\subset U$ and 
such that 
for each $(i,j)$ with $i\neq j$, $f_{i}^{-1}(U)\cap 
f_{j}^{-1}(U)=\emptyset $, 
then $\HD(J(G))=\d (f)$, 
where $f=(f_{1},\ldots ,f_{s}).$ 
As a direct consequence of Theorem A, 
the statement of Theorem B holds. 
\qed 

\

\section{Examples}
\label{Examples}
Throughout this section, 
we provide 
an extensive collection of classes of examples of analytic family 
of semigroups satisfying all the hypothesis of Theorem A and Theorem B 
and we analyze 
in detail the corresponding Bowen parameter or Hausdorff dimension 
function.    

First, we give some examples of analytic families of 
expanding rational semigroups 
satisfying the open set condition. 
\begin{ex}
Let $d_{1},d_{2}\in \Bbb{N}, \geq 2$ with 
$(d_{1},d_{2})\neq (2,2).$ Let $a\in \Bbb{C}$ with 
$0<|a|<1.$ Then, setting $g=(z^{d_{1}},az^{d_{2}}),$ 
we see that the semigroup $G=\langle g_{1},g_{2}\rangle $ is hyperbolic, and hence
 expanding. 
There exists also an open neighborhood $U$ of 
$\{ z\in \Bbb{C}: 1\leq |z|\leq 1/|a|^{\frac{1}{d_{2}-1}}\} $ such that
$g_{1}^{-1}(\overline{U})\cup g_{2}^{-1}(\overline{U})\subset U \mbox{ and } 
g_{1}^{-1}\overline{U})\cap g_{2}^{-1}(\overline{U})=\emptyset . 
$
Hence, small perturbation of $g$ satisfies the same property. Therefore,  
if we take a small neighborhood of $V$ of $g$ in 
$(\mbox{{\em Rat}})^{2}$, then setting 
$G_{f}:= \langle f_{1},f_{2}\rangle $ for each $f\in V$, 
we have that $\{G _{f}\} _{f\in V}$ is an analytic family of 
expanding rational semigroups and for each $f\in V$, the semigroup 
$G_{f}$ satisfies the strongly separating open set condition with 
the set $U.$  
By Remark~\ref{oscrem}, for each $f\in V$, 
$\delta (f)=\HD(J(G_{f}))<2.$    
\end{ex}
\begin{prop}(See \cite{sumid1, sumi07}) 
Let $f_{1}$ be a hyperbolic polynomial 
with $\deg (f_{1})\geq 2$ such that 
$J(f_{1})$ is connected. 
Let $K(f_{1})$ be the filled-in Julia set of $f_{1}$ and 
let $b\in \mbox{{\em int}} K(f_{1})$ be a point. 
Let $d$ be a positive integer such that 
$d\geq 2.$ Suppose that $(\deg (f_{1}),d)\neq (2,2).$ 
Then, there exists a number $c>0$ such that 
for each $\l \in \{ \l\in \Bbb{C}: 0<|\l |<c\} $, 
setting $f_{\l }=(f_{\l ,1},f_{\l ,2})=
(f_{1},\l (z-b)^{d}+b ),$ we have that 
$f_{\l }\in \Exp(2)$,  
$f_{\l }$ satisfies the separating open set condition with 
an open set $U_{\l }$, 
$J(\langle f_{\l ,1},f_{\l ,2}\rangle )$ is porous,  
$\HD(J(\langle f_{\l 1.},f_{\l, 2}\rangle ))=\delta 
(f_{\l })<2$, and $P(\langle f_{\l ,1}, f_{\l ,2}\rangle )
\setminus \{ \infty \} $ is bounded in $\Bbb{C}.$  
\end{prop}
{\sl Proof.} 
We will follow the argument in \cite{sumid1, sumi07}. 
Conjugating $f_{1}$ by a M\"{o}bius transformation, 
we may assume that $b=0$ and the coefficient 
of the highest degree term of $f_{1}$ is equal to $1.$ 

 For each $r>0$, we denote by $D(0,r)$ the Euclidean disc with 
 radius $r$ and center $0.$ Let $r>0$ be a 
 number such that $\overline{D(0,r)}\subset \mbox{int} K(f_{1}).$ 
We set $d_{1}:=\deg (f_{1}).$ 
 Let $\alpha >0$ be a number. 
Since $d\geq 2$ and $(d,d_{1})\neq (2,2)$, 
it is easy to see that 
$(\frac{r}{\alpha })^{\frac{1}{d}}>
2\left(2(\frac{1}{\alpha })
^{\frac{1}{d-1}}\right)^{\frac{1}{d_{1}}}
$ if and only if 
\begin{equation}
\label{Contproppfeq1}
\log \alpha <
\frac{d(d-1)d_{1}}{d+d_{1}-d_{1}d}
( \log 2-\frac{1}{d_{1}}\log \frac{1}{2}-\frac{1}{d}\log r) .
\end{equation} 
We set 
\begin{equation}
\label{Contproppfeq2}
c_{0}:=\exp \left(\frac{d(d-1)d_{1}}{d+d_{1}-d_{1}d}
( \log 2-\frac{1}{d_{1}}\log \frac{1}{2}-\frac{1}{d}\log r) \right)
\in (0,\infty ).
\end{equation}
 
Let $0<c<c_{0}$ be a small number and let $\l \in \Bbb{C} $ 
be a number with $0<|\l |<c.$ 
Put $f_{\l ,2}(z)=\l z^{d}.$ 
Then, we obtain $K(f_{\l ,2})=\{ z\in \Bbb{C} \mid 
|z|\leq (\frac{1}{|\l |})^{\frac{1}{d-1}}\} $ and 
$$f_{\l,2}^{-1}(\{ z\in \Bbb{C} \mid  |z|=r\} )=
\{ z\in \Bbb{C} \mid |z|=(\frac{r}{|\l |})^{\frac{1}{d}}\} .$$
Let 
$D_{\l }:=\overline{D(0,2(\frac{1}{|\l |})^{\frac{1}{d-1}})}.$ 
Since $f_{1}(z)=z^{d_{1}}(1+o(1))\ (z\rightarrow \infty )$,  
it follows that if $c$ is small enough, then 
for any $\l \in \Bbb{C} $ with $0<|\l |<c$, 
$$f_{1}^{-1}(D_{\l })\subset 
\left\{ z\in \Bbb{C} \mid 
|z|\leq 2\left( 2(\frac{1}{|\l |})^{\frac{1}{d-1}}\right) 
^{\frac{1}{d_{1}}}\right\} .$$  
This implies that  
\begin{equation}
\label{Contproppfeq3}
f_{1}^{-1}(D_{\l })\subset f_{\l ,2}^{-1}(\{ z\in \Bbb{C} \mid |z|<r\} ).
\end{equation} 
Hence, 
setting $U_{\l }:=\mbox{{\em int}} K(f_{\l ,2})\setminus K(f_{1})$, 
$f_{\l }=(f_{1},f_{\l ,2})$ satisfies the 
separating open set condition with $U_{\l }.$  
Therefore, setting $G_{\l }:=\langle f_{1}, f_{\l, 2}\rangle $, 
we have $J(G_{\l })\subset \overline{U_{\l }}
\subset K(f_{\l ,2})\setminus \mbox{{\em int}} K(f_{1}). $ 
In particular, {\em int}$K(f_{1})\cup (\oc \setminus K(f_{\l ,2}))\subset 
F(G_{\l }).$ 
Furthermore, (\ref{Contproppfeq3}) implies that 
$f_{\l ,2}(K(f_{1}))\subset \mbox{{\em int}} K(f_{1}).$ 
Thus, 
we have $P(G_{\l })\setminus \{ \infty \} 
= \overline{\bigcup _{g\in G_{\l }\cup \{ Id\} }
g(CV^{\ast }(f_{1})\cup CV^{\ast }(f_{\l, 2}))}
\subset 
\mbox{{\em int}} K(f_{1})\subset F(G_{\l })$, 
where $CV^{\ast }(\cdot )$ denotes the set of 
all critical values in $\Bbb{C}.$ Hence, 
$G_{\l }$ is expanding and 
$P(G_{\l })\setminus \{ \infty \} $ is bounded in $\Bbb{C}.$ 
By Theorem B and Remark~\ref{oscrem}, we obtain that 
for each $\l $ with $0<|\l |<c$, 
$J(G_{\l })$ is porous and  
$\HD(J(G_{\l }))=\delta (f_{\l })<2.$ We are done.  
\qed 
\begin{ex}[\cite{sumi07}]
Let $h_{1}(z)=z^{2}/4,h_{2}(z)=z^{2}-1,f_{1}:=h_{1}^{2},f_{2}:=h_{2}^{2},$ 
and 
$f:=(f_{1},f_{2}).$ Let $G=\langle f_{1},f_{2}\rangle .$ 
Then it is easy to see that 
$f_{1}(K(f_{2}))\subset \mbox{{\em int}}(K(f_{2}))$ and 
$P(G)\setminus \{ \infty \} \subset \mbox{{\em int}}(K(f_{2})).$ 
Hence, we have $P(G)\subset F(G)$, which implies that $f\in \Exp(2).$ 
Moreover, it is easy to see that $f$ satisfies the strongly separating 
open set condition with an open set $U$ (letting $U$ be an open 
neighborhood of $K(f_{1})\setminus \mbox{{\em int}}(K(f_{2}))$).  
Thus there exists an open neighborhood 
$V$ of $f$ in $(\mbox{{\em Rat}})^{2}$ 
 such that for each $g\in V$, we have that 
$g\in \Exp(2)$, $g$ satisfies the strongly separating 
open set condition with $U$, 
$P(\langle g_{1},g_{2}\rangle )\setminus \{ \infty \} $ is bounded in 
$\C $, and 
$\HD(J(\langle g_{1},g_{2}\rangle ))=\delta (g)<2.$  
See Figure~\ref{fig:dcjulia} for the Julia set of 
$\langle f_{1},f_{2}\rangle .$ 
\begin{figure}[htbp]
\caption{The Julia set of 
$\langle f_{1},f_{2}\rangle$, 
where $h_{1}(z):=z^{2}/4,\ h_{2}(z):=z^{2}-1,\ f_{1}:=h_{1}^{2},\ 
f_{2}:=h_{2}^{2}.$}
\includegraphics[width=5cm,width=5cm]{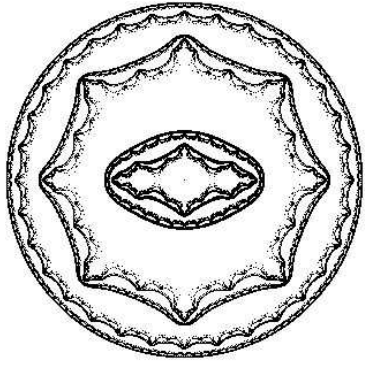}
\label{fig:dcjulia}
\end{figure}

\end{ex}
\begin{ex}
For each $j=1,2,$ 
let $\gamma _{j}$ be a hyperbolic polynomial  
such that $\deg (\gamma _{j})\geq 2$ and 
$J(\gamma _{j})$ is connected.  
Suppose that $K(\gamma _{1})\cap K(\gamma _{2})=\emptyset $, 
where $K(\cdot )$ denotes the filled-in Julia set. 
Let $R>0$ be a large number such that 
$B(0,R)\supset K(\gamma _{1})\cup K(\gamma _{2}).$ 
Then, there exists a large positive integer $n$ such that 
with $U:= B(0,R)$, 
\begin{equation}
\label{oscg}
\gamma _{1}^{-n}(\overline{U})\cup 
\gamma _{2}^{-n}(\overline{U})
\subset U \mbox{ and } 
\gamma _{1}^{-n}(\overline{U})\cap 
\gamma _{2}^{-n}(\overline{U})=\emptyset .
\end{equation} 
Thus, setting $g=(g_{1},g_{2}):=(\gamma _{1}^{n}, \gamma _{2}^{n})$,
 there exists an open neighborhood 
 $V$ of $g$ in $(\mbox{{\em Rat}})^{2}$ such that 
each $f=(f_{1},f_{2})\in V$ satisfies the strongly separating open set condition with 
$U$. Since 
$K(g_{j})\subset U$ for each $j=1,2$, 
(\ref{oscg}) implies that 
setting $W:= \mbox{{\em int}}K(g_{1})\cup \mbox{{\em int}}K(g_{2})
\cup (\oc \setminus \overline{U})$, 
we have that for each $j=1,2$, $g_{j}(W)\subset W.$ 
Hence, $W\subset F(\langle g_{1},g_{2}\rangle ).$ 
Combining this, $CV^{\ast }(g_{j})\subset K(g_{j})\subset U$
where $CV^{\ast }(g_{j} )$ denotes the set of 
critical values of $g_{j}$ in $\Bbb{C}$, 
 and 
(\ref{oscg}), we obtain that 
$P(\langle g_{1},g_{2}\rangle )\subset 
F(\langle g_{1},g_{2}\rangle ).$ 
Therefore, $g=(g_{1},g_{2})\in \Exp(2).$ 
Thus, if we take the above $V$ small enough, 
it follows that 
for each $f=(f_{1},f_{2})\in V$, 
we have that $f\in \Exp(2)$ and $f$ satisfies the 
strongly separating open set condition. In particular, for each $f\in V$, 
$\HD(J(\langle f_{1},f_{2}\rangle ))=\delta (f)<2$ and 
$J(\langle f_{1},f_{2}\rangle )$ is porous.  
\end{ex}

Now we describe an example of analytic family 
$\{ G_{\l }\} _{\l \in \Lambda }$ 
of expanding rational semigroups, 
where $G_{\l }=\langle f_{\l ,1},\ldots ,f_{\l ,s}\rangle $, 
$f_{\l }=(f_{\l ,1}, \ldots ,f_{\l ,s})$, 
such that each $G_{\l }$ satisfies the open 
set condition with $U_{\l }$ 
but does not satisfy the separating open set condition
with any open subset of $\oc $, and 
such that for each $\l $, $\HD(J(G_{\l }))=\delta (f_{\l })<2.$  
\begin{ex}(See \cite[Example 6.2]{sumi2}) 
Let $p_{1},p_{2},$ and $p_{3}\in \Bbb{C}$ be mutually distinct 
points that form an equilateral triangle. 
Let $U$ be the interior part of the triangle. 
Let $\gamma _{j}(z)=2(z-p_{j})+p_{j}$, for each $j=1,2,3.$ 
Then 
$J(\langle \gamma _{1},\gamma _{2},\gamma _{3}\rangle )$ 
is the Sierpi\'{n}ski gasket, which is connected. 
Hence $(\gamma _{1}, \gamma _{2},\gamma _{3})$ 
satisfies the open set condition with $U$ but 
fails to  
satisfy the separating open set condition with any open 
subset of $\oc .$ 
 Let $x$ be the barycenter of the equilateral triangle $p_{1}p_{2}p_{3}$ 
and let $r>0$ be a small number such that 
$D(x,r)\subset U\setminus \cup _{j=1}^{3}\gamma _{j}^{-1}(\overline{U})$, 
where $D(x,r)$ denotes the Euclidean disc with center $x$ and radius $r.$ 
Let $\gamma _{4}$ be a polynomial such that 
$J(\gamma _{4})=\partial D(x,r).$ 
Let $u$ be a large positive integer such that 
$\gamma _{4}^{-u}(\overline{U})\subset 
U\setminus \cup _{j=1}^{3}\gamma _{j}^{-1}(\overline{U}).$ 
Set $\alpha := \gamma _{4}^{u}.$ Then, 
there exists a neighborhood $V$ of $\alpha $ in $\mbox{{\em Rat}}$ 
such that for each $\beta \in V$, 
$f=(\gamma _{1},\gamma _{2}, \gamma _{3},\beta )$ satisfies
 the open set condition with $U.$ 
Let $G_{\beta }=\langle \gamma _{1},\gamma _{2},
 \gamma _{3}, \beta \rangle $, 
 for each $\beta \in V.$  
If we take $u$ large enough, 
 then  
 we may assume that $P(G_{\alpha })\subset F(G_{\alpha }).$ Therefore, 
 $G_{\alpha }$ is expanding. Thus, if we take the above $V$  small enough, 
 it follows that 
 for each $\beta \in V$, 
we have  that $G_{\beta }$ is expanding, 
$G_{\beta }$ satisfies the open set condition with $U$, 
and $G_{\beta }$ fails to satisfy the separating open set 
condition with any open subset of $\oc .$ 
also note that for each $\beta \in V$, 
$\cup _{j=1}^{3}\gamma _{j}^{-1}(\overline{U})
\cup \beta ^{-1}(\overline{U})$ is a proper subset of 
$\overline{U}.$ Combining it with \cite[Lemma 2.4]{hiroki1} and 
\cite[Theorem 1.25]{sumi06}, 
we obtain that for each $\beta \in V$, 
$J(G_{\beta })$ is porous and $\HD (J(G_{\beta }))<2.$ Thus, 
for each $\beta \in V$, we have 
$$\HD(J(G_{\beta }))=\delta 
(\gamma _{1},\gamma _{2},\gamma _{3}, \beta)<2.$$
See Figure~\ref{fig:sierpjulia} for the Julia set of $G_{\alpha }.$     
\begin{figure}[htbp]
\caption{The Julia set of 
$G_{\alpha }=\langle \gamma _{1},\gamma _{2},\gamma _{3}, 
\alpha \rangle .$}
\includegraphics[width=6cm,width=6cm]{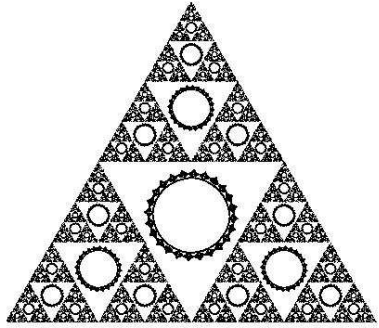}
\label{fig:sierpjulia}
\end{figure}
\end{ex}
\begin{rem}
In the sequel \cite{sumiprepare2} (announced in \cite{sumikokyuroku}), 
we will see that there are plenty of parameters 
$f=(f_{1},f_{2})\in \Exp(2)$ such that $f$ satisfies 
all of the following conditions: 
(1)$f_{1}$ and $f_{2}$ are 
polynomials of degree greater than or equal to two, (2)$f$ satisfies 
the open set condition, (3)$P(\langle f_{1},f_{2}\rangle )\setminus 
\{ \infty \} $ is bounded in $\C $, 
(4)$\HD(J(\langle f_{1},f_{2}\rangle ))<2$, 
and (5)$J(\langle f_{1},f_{2}\rangle )$ is connected. 
\end{rem}
We give an example of analytic family 
$\{ G_{\l }\} _{\l \in \Lambda }
$ 
of expanding rational semigroups 
satisfying the open set condition and 
$\delta (f_{\l })=\HD(J(G_{\l }))=2$, 
where $G_{\l }=\langle f_{\l ,1},\ldots ,f_{\l ,s}\rangle , f_{\l }:=
(f_{\l ,1},\ldots ,f_{\l ,s})$, 
$\Lambda =\{ \l\in \Bbb{C}: 0<|\l|<1\} .$  
\begin{ex}
Let $\l\in \Bbb{C}$ with $0<|\l|<1$ and 
let $f_{\l }=(f_{\l ,1},f_{\l ,2})=(z^{2},\l z^{2})\in (\mbox{{\em Rat}})^{2}.$ 
Let $G_{\l }=\langle f_{\l, 1},f_{\l, 2}\rangle .$ Then, 
$P(G_{\l })=\{ 0,\infty \} \subset F(G_{\l })$ and therefore 
$G_{\l }$ is expanding. Let 
$A_{\l }=\{ z\in \Bbb{C}: 1< |z|< 1/|\l |\} .$ 
Then, we have 
$f_{\l, 1}^{-1}(A_{\l })\cup f_{\l, 2}^{-1}(A_{\l })\subset A_{\l }$ and 
$f_{\l ,1}^{-1}(A_{\l })\cap f_{\l,2}^{-1}(A_{\l })=\emptyset .$ Hence, 
$G_{\l }$ satisfies the open set condition. Therefore, 
by \cite{sumi2}, we have $\delta (f_{\l })=\HD(J(G_{\l })).$ 
Since the point $1$ belongs to the Julia set of $f_{1}=f_{\l ,1}$, 
we have $1\in J(G_{\l }).$ Moreover, we easily obtain that 
$\cup _{g\in G_{\l }}g^{-1}(1) $ is dense in $\overline{A_{\l }}.$ 
Hence, by \cite[Lemma 3.2]{HM}, it follows that 
$J(G_{\l })=\overline{A_{\l }}.$ Thus, 
$\delta (f_{\l })=\HD(J(G_{\l }))=2.$  

\end{ex}
We give another example of analytic family of expanding rational 
semigroups such that $\HD(J(\langle f_{\l, 1},\ldots 
,f_{\l s}\rangle ))\leq 
\delta (f_{\l })<2.$ 
\begin{ex}
\label{d1dsex1}
Let $d_{1},d_{2},\ldots ,d_{s}\in \Bbb{N}$ such that 
$d_{j}\geq 2$ for each $j=1,\ldots ,s.$ 
Let $g=(g_{1},\ldots ,g_{s})=(z^{d_{1}},\ldots ,z^{d_{s}})\in 
(\mbox{{\em Rat}})^{s}.$ 
Let $L_{t}:C(J(\tilde{g}))\rightarrow C(J(\tilde{g}))$ be the 
Perron Frobenius operator defined by the formula
$$L_{t}\varphi (z)=\sum _{\tilde{g}(y)=z}|\tilde{g}'(y)|^{-t}\varphi (y).$$ 
Then, we have $L_{t}1\equiv (\sum _{j=1}^{s}\frac{1}{d_{j}^{t-1}})1$, 
where $1$ denotes the constant function taking the value $1.$ 
Hence, setting $\beta (t)=\sum _{j=1}^{s}\frac{1}{d_{j}^{t-1}},$ 
we have 
\begin{equation}
\label{betap}
\beta (t)=e^{P(t,g)}.
\end{equation}
 Thus, 
$\beta (\delta (g))=1.$ 
We now assume that 
$\sum _{j=1}^{s}\frac{1}{d_{j}}<1.$ 
Then, since the function $t\mapsto \beta (t)$ is strictly decreasing, 
we obtain that $\delta (g)<2.$ 
Since $f\mapsto \delta (f)$ is continuous around $g$, 
it follows from \cite{sumi2} that 
there exists an open neighborhood $U$ of $g$ 
in $(\mbox{{\em Rat}})^{s}$ and an $\epsilon >0$ 
such that 
for each $f\in U$, 
$\HD(\langle f_{1},\ldots ,f_{s}\rangle )
\leq \delta (f)\leq 2-\epsilon .$ In particular, for each 
$f\in U$, 
{\em int}$J(\langle f_{1},\ldots ,f_{s}\rangle ))=\emptyset .$ 
By Remark~\ref{strem} and Remark~\ref{dimjexppaperrem}, 
for almost every $f\in U$ with respect to the Lebesgue measure, 
$s_{0}(\langle f_{1},\ldots ,f_{s}\rangle )=
t_{0}(f)=\delta (f)\leq 2-\epsilon <2.$   

\end{ex}
Let us now provide several sufficient conditions for $f$ to satisfy $\delta (f)>2.$ 
\begin{ex}
Using the same notation as that in Example~\ref{d1dsex1}, 
suppose there exists an integer $m$ such that 
$d_{1}=\cdots =d_{s}=m.$ 
Then, by (\ref{betap}), 
we obtain $\delta (g)=1+\frac{\log s}{\log m}.$ 
Since the function $f\mapsto \delta (f)$ is continuous and plurisubharmonic 
around $g$, 
it follows that for each open neighborhood $U$ of 
$g$ in $(\mbox{{\em Rat}})^{s}$, there exists 
a non-empty open subset $V$ of $
U\setminus \{ g\} $ such that for each $f\in V $,  
\begin{equation}
\label{d1dsex2eq1}
\delta (f)\geq 1+\frac{\log s}{\log m}.
\end{equation} 

 We now assume that $s>m.$ Then, 
 from the equality
 $\delta (g)=1+\frac{\log s}{\log m}$  and the continuity of 
 $f\mapsto \delta (f)$ around $g$, 
 it follows that for each $\epsilon $ with 
 $0<\epsilon <\frac{\log s}{\log m}-1$, 
there exists an open neighborhood $W$ of 
$g$ in $(\mbox{{\em Rat}})^{s}$ such that 
for each $f\in W$, 
$\delta (f)\geq 1+\frac{\log s}{\log m}-\epsilon >2.$ 
In particular, 
for each $f\in W$, $f$ does not satisfy 
the open set condition. Moreover, 
by Remark~\ref{strem} and Remark~\ref{dimjexppaperrem}, 
for almost every $f\in W$ with respect to the Lebesgue measure, 
we have $s_{0}(\langle f_{1},\ldots , f_{s}\rangle )=
t_{0}(f)=\delta (f)\geq 1+\frac{\log s}{\log m}-\epsilon >2.$   
Note that for a fixed $m$, 
$1+\frac{\log s}{\log m}\rightarrow \infty $ as 
$s\rightarrow \infty .$ Thus, 
the functions $f\mapsto \delta (f)$, 
$f\mapsto t_{0}(f), $ and 
$f\mapsto s_{0}(\langle f_{1},\ldots ,f_{s}\rangle )$, 
where $f\in \Exp(s)$, are unbounded, 
 if $s$ runs over all 
positive integers.  

\end{ex}
\begin{prop}
\label{deltageq2prop}
Let $g=(g_{1},\ldots ,g_{s})\in \Exp(s)$ and let 
$G=\langle g_{1},\ldots ,g_{s}\rangle .$ Let 
$m_{2}$ be the $2$-dimensional Lebesgue measure. 
Suppose that 
there exists a couple $(i,j)$ with $i\neq j$ such that 
$m_{2}(g_{i}^{-1}(J(G))\cap 
g_{j}^{-1}(J(G)))>0.$  
Then, $\delta (g)>2.$ 
In addition, for each $0<\epsilon <\delta(f)-2$, 
there exists an open neighborhood $U$ of $g$ in $(\mbox{{\em Rat}})^{s}$ such that 
for each $f\in U$, $\delta (f)\geq \delta (g)-\epsilon >2.$ 
\end{prop}
{\sl Proof.}
By the assumption, \cite{sumi2} implies that 
$\delta (g)\geq \HD(J(G))=2.$ Suppose $\delta (g)=2.$ Then, 
by \cite[Proposition 4.13]{sumi2}, we obtain a contradiction. 
Thus, $\delta (g)>2.$ Since the function $f\mapsto \delta (f)$ is 
continuous around $g$, the rest of the statement of the proposition 
holds.   
We are done. 
\qed 

\begin{ex}
Let $g=(g_{1},g_{2},g_{3})=(z^{2},z^{2}/4,z^{2}/3)\in (\mbox{{\em Rat}})^{3}$ 
and let $G:=\langle g_{1},g_{2},g_{3}\rangle .$  
Then, $P(G)=\{ 0,\infty \} \subset F(G)$ and so 
$G$ is expanding and $J(G)=\{ z\in \Bbb{C}: 1\leq |z|\leq 4\} .$ 
By \cite[Example 4.14]{sumi2}, 
(or by Proposition~\ref{deltageq2prop}), 
we have that $\delta (g)>2.$ 
Since the function $f\mapsto \delta (f)$ is continuous around $g$, 
it follows that 
for each $\epsilon $ with $0<\epsilon <\delta (g)-2$, 
there exists an open neighborhood $U$ of $g$ in $(\mbox{{\em Rat}})^{3}$ such that
 for each $f\in U$, 
 $\delta (f)\geq \delta (g)-\epsilon >2.$ 
Moreover, by Remark~\ref{strem} and Remark~\ref{dimjexppaperrem}, 
for almost every $f\in U$ with respect to the Lebesgue measure, 
we have $s_{0}(\langle f_{1},f_{2},f_{3}\rangle )
=t_{0}(f)=\delta (f)\geq \delta (g)-\epsilon >2.$ 
\end{ex}

\end{document}